\documentclass[a4paper,twoside,11pt]{amsart}

\usepackage{my_package}

\title{A note on ideals in derived geometries}
\date{\today}

\author{Zachary Gardner}
\address{Boston College\\Department of 
	Mathematics\\ Chestnut Hill MA 02467 \\USA}
\email{zachary.gardner@bc.edu}

\author{Jeroen Hekking} 
\address{University of Regensburg\\Department of Mathematics\\Regensburg 93053\\Germany}
\email{jeroen.hekking@ur.de}

\begin{document}
	\maketitle

\begin{abstract}
    We develop the basic theory of derived quasi-coherent ideals for stacks relative to a given derived algebraic context. We compare different notions of adic completeness with respect to derived ideals, define and compare formal spectra and formal completions along closed immersions, and connect the theory of derived ideals to that of derived extended Rees algebras. A first application is the construction of derived scheme-theoretic images in full generality. We further show that the deformation space of any nonconnectively affine morphism of derived stacks is nonconnectively affine over the base. We close with a first exploration of transmutation cohomology and filtrations thereof in this more general context.
\end{abstract}
\setcounter{tocdepth}{1}
\tableofcontents

 \section*{Introduction}
Let $A$ be an animated ring, let $I \to A$ be a generalized Cartier divisor in the sense of \cite{KhanVirtual}, and let $M$ be a derived $A$-module in the unbounded derived $\infty$-category of $A$. The start of this project has simply been:
\begin{Quest**}
    What does it mean for $M$ to be $I$-complete?
\end{Quest**}
This basic question has relevance in the current study of derived prismatic cohomology, see \cite{MaoPrismaticlogarithm, HolemanDerivedDelta, BhattLuriePrismatization, BhattAbsolute}, and answers can be found for example in \cite{LurieHA,DAGXII,DwyerComplete,GreenleesDerived}. In this note we compare the following two answers, also in a more general setting:

\begin{Ansr**} 
    Say that $M$ is \emph{derived $I$-complete} if it is local with respect to the  $I$-equivalences. Here, \emph{$I$-equivalences} are those maps $N' \to N$ such that $N' \otimes_A I \to N \otimes_A I$ becomes invertible.
\end{Ansr**}

\begin{Ansr**} 
\label{Anr_B}
    Say that $M$ is \emph{derived $I$-adically complete} if the natural map $M \to M^\wedge_I$ is invertible. Here, $M^\wedge_I$ is the \emph{derived $I$-adic completion} of $M$ which is a derived analogue of the formula $\lim M/I^nM$ and is computed via the derived extended Rees algebra. 
\end{Ansr**}

The different notions of $I$-completeness proposed in these answers are completely coherent for any derived ideal pair $(A,I)$, meaning that $I \to A$ is a derived Smith ideal in the sense of \cite{MagidsonDividedpowers}.  In these notes, we  generalize this to derived ideal pairs $(A,I)$ relative to a given---and fixed in this section---derived algebraic context $\CCC$ as in \cite[Def.~4.2.1]{RaksitDerived}, see \S\ref{Subsec:Derived_ideal_pairs}. This uses the full force of the theory of derived algebraic contexts: a derived ideal pair over $\CCC$ is nothing but a derived ring over an appropriate derived algebraic context supported on the arrow category $\Fun(\Delta^1,\CCC)$. We globalize this to \emph{derived quasi-coherent ideals} $\II$ on derived $\CCC_{\geq 0}$-stacks $X$, see Definition~\ref{Def:I_adic_complete} and Definition~\ref{Def:I_complete}. We then show:

\begin{Thm**}[\textbf{Completeness}]
\label{ThmB}
    Let $X$ be a $\CCC_{\geq 0}$-stack, $\II$ a quasi-coherent ideal on $X$, and $\MM$ a quasi-coherent $\OO_X$-module. If $\MM$ is derived $\II$-adically complete then it is derived $\II$-complete (Corollary~\ref{Cor:derived_complete_implies_adically}). In the case that $\II$ is locally finitely generated, then also the converse is true (Theorem~\ref{Thm:derived_I_complete_is_adically}). 
\end{Thm**}

Here, we use the theory of \emph{$\CCC_{\geq 0}$-stacks} and their derived quasi-coherent $\OO_X$-modules and $\OO_X$-algebras as developed in \cite{OrenJeroenblowups}, see \S \ref{subsec:geometric_cntxt}. To make sense of $\MM^{\wedge}_\II \triangleq \lim \MM/\II^n\MM$ for $\MM \in \QCoh(X)$ in Answer~\ref{Anr_B} we make use of the \emph{derived $\II$-adic filtration} on $\MM$ in Definition~\ref{Def:I_adic_filtration}. This is induced by tensoring with the derived extended Rees algebra $\RR_{V(\II)/X}^\ext$ associated to the derived vanshing locus $V(\II)$ attached to $\II$, where $V(\II)$ is (the restriction of) the relative nonconnective spectrum of $\OO_X/\II$ over $X$.\footnote{See Remark~\ref{Rem:affine_stacks} for an overview of some different notions of affineness.}

In algebraic terms and affine-locally, the \emph{derived vanishing locus} is one side of a beautiful generalization of the classical correspondence between surjections $A \to B$ of discrete $A$-algebras on one hand, and ideals $I$ of $A$ on the other. In the derived world one has access to not only animated structures but also to nonconnective homotopical data. Translating this to derived rings and globalizing yields:
\begin{Thm**}[\textbf{Affines-as-ideals}]
    Let $X$ be a $\CCC_{\geq 0}$-stack. Then there is an equivalence between derived quasi-coherent ideals $\II$ of $X$ and derived quasi-coherent $\OO_X$-algebras $\BB$, induced by sending $\II \to \OO_X$ to its cofiber and conversely by sending $\OO_X \to \BB$ to its fiber (Lemma~\ref{Lem:Idea_QAlg}). Restricting to connectives and passing to associated relative spectra induces an equivalence between closed immersions $Z \to X$ and connective quasi-coherent ideals $\II$ of $X$ (Corollary~\ref{Cor:V_Cld_equiv}).
\end{Thm**}
The essential ingredient in this theorem is the affine statement, which says that the category of derived ideal pairs $(A,I)$ is equivalent to the category of arrows $A \to B$ in $\DAlg$.\footnote{The fact that this `only' translates to a correspondence between closed immersions and connective quasi-coherent ideals on the global level is because we are working with $\CCC_{\geq 0}$-stacks. This should disappear when going to the nonconnective theory of $\CCC$-stacks, although $t$-structures---and by extension closed immersions---can be much more subtle here, see \cite{MathewAffine}. With an eye towards applications in `animated' algebraic geometry (that is, the connective part), we restrict ourselves to $\CCC_{\geq 0}$-stacks.}  In fact, this is a reflection of a deeper statement on the level of underlying derived algebraic contexts, see \S \ref{Subsec:Derived_ideal_pairs}. For $\CCC = \Mod_\Z$ this has recently also been observed in \cite[\S 3.1]{MagidsonDividedpowers}, and the main ingredients are  in \cite[\S 3.3]{MaoRevisiting}. 

\textbf{A first application} is the construction of the \emph{derived scheme-theoretic image} of an arbitrary morphism $f \colon W \to X$ of $\CCC_{\geq 0}$-stacks. This yields a factorization of $f$ through a closed immersion $i \colon Z \to X$ which is initial in the category of factorizations of $f$ of the form
\[ W \xrightarrow{f'} Z' \xrightarrow{i'} X, \]
where $i'$ is a closed immersion (Proposition~\ref{Prop:scheme_theoretic_image}). The construction mimics the classical picture: one defines $Z$ as the relative spectrum over $X$ of the image $\varphi \OO_X$ of the canonical map $\varphi \colon \OO_X \to f_* \OO_W$, where $\varphi \OO_X$ is induced using connective covers of derived ideal pairs. When $f$ is affine, then $f_* \OO_W$ is connective, in which case the induced map $\varphi \OO_X \to f_* \OO_X$ is injective on $\pi_0$ and invertible on $\pi_n$ for $n>0$ (Lemma~\ref{Lem:image_ring_map}). In general $f_*\OO_X$ need not be connective, hence the need for nonconnective derived geometry in this construction. In the context $\CCC = \Mod_\Q$ one can work with $\CAlg(\CCC)$ instead, as done in \cite{GaitsgoryStudy}(Example~\ref{Ex:GR_image}).  

\textbf{A second application} is the theory of formal spectra in geometry relative to $\CCC$. For a quasi-coherent ideal $\II$ on a $\CCC_{\geq 0}$-stack $X$, we define the \emph{formal spectrum} $\Spf(\OO_X,\II)$ as the colimit over all infinitesimal neighborhoods $V(\II^n) \coloneqq \widetilde{\Spec}_X (\OO_X/\II^n)$, meanwhile the \emph{formal completion} $X^\wedge_{\lvert V(\II) \rvert}$ along $V(\II) \to X$ is defined as the largest subobject of $X$ which has empty intersection with $X \setminus V(\II)$.\footnote{Following \cite{MathewAffine}, we write $\widetilde{\Spec}$ for the functor from derived $\CCC$-algebras to $\CCC_{\geq 0}$-stacks induced by restricted Yoneda. This is not fully faithful in general: see Warning~\ref{Warn:nc_affine} and references there.} We then show:
\begin{Thm**}[\textbf{Formal spectrum is formal completion}]
    Suppose that the quasi-coherent ideal $\II$ on $X$ is locally finitely generated. Then the canonical map
    \[ \Spf(\OO_X,\II) \to X^\wedge_{\lvert V(\II) \rvert} \]
    is invertible.
\end{Thm**}

The present notes also finish the story of extended Rees algebras in derived algebraic geometry which started with the construction in the case of animated rings $A \to B$ which are surjective on $\pi_0$ over $\CCC = \Mod_\Z$ in \cite{HekkingGraded}, and which was generalized to arbitrary affine morphisms $X \to Y$ of $\CCC_{\geq 0}$-stacks in \cite{OrenJeroenblowups} for $\CCC$ arbitrary. In the latter case, the derived extended Rees algebra $\RR_{X/Y}^\ext$ is in general not connective, but the input is still required to be so. We generalize this to also allow morphisms of derived $\CCC$-algebras $A \to B$ which are nonconnective (Definition~\ref{Def:Rees}), and use the global version of this construction to show:
\begin{Thm**}[\textbf{Affines are closed under deformation spaces}]
\label{ThmF}
    Let $\BB$ be a quasi-coherent algebra on $X$ and $W \coloneqq \Spec_X(\BB)$ its relative (nonconnective or almost affine) spectrum. Then the  deformation space $D_{W/X}$ is (nonconnectively or almost) affine over $X \times \A^1$, and in fact $D_{W/X} \simeq \Spec_X(\RR_{W/X}^\ext)$ (Corollary~\ref{Cor:affine_def}, Remark~\ref{Rem:nc_affine_D}).\footnote{Here, the \emph{($\G_m$-equivariant) deformation space} is as in \cite[Thm.~A]{Weil}, and our result generalizes \cite[Thm.~4.6]{OrenJeroenblowups}.}
\end{Thm**}
Let $A \to B$ be a map of discrete rings over $\Mod_\Z$, with ideal $I$. Recall that $\pi_1\LL_{B/A} \simeq I/I^2$ \cite[Cor.~3.14]{QuillenHomology}. Using derived quasi-coherent ideals and the deformation to the normal bundle, we show:
\begin{Thm**}[\textbf{Cotangent complex as derived conormal sheaf}]
\label{ThmQ}
    Let $Z \to X$ be a closed immersion of $\CCC_{\geq 0}$-stacks, with corresponding ideal $\II$. Then there is a canonical equavalence
    \[ \LL_{Z/X}[-1] \simeq \II/\II^2, \]
    where $\II^2$ is the ideal associated to the derived first order infinitesimal neighborhood $Z^{(1)} \to X$ of $Z \to X$  (Remark~\ref{Rem:inf_nbhd}, Corollary~\ref{Cor:cotangent_Quillen}).
\end{Thm**} 
In fact, the proof of Theorem~\ref{ThmQ} is a straightforward application of the deformation to the normal bundle to the closed immersion $Z \to X$, which implies that
\[ \LSym_{\OO_Z}(\LL_{Z/X}[-1](-1)) \simeq \bigoplus_{n \in \N} \II^n/\II^{n+1}(-n), \]
where $\II^{k+1}$ is the quasi-coherent ideal associated to the derived $k$-th order infinitesimal neighborhood $Z^{(k)} \to X$. It follows that the derived normal bundle can also be understood as a derived analogue of the normal cone (using the convention followed in \cite[\href{https://stacks.math.columbia.edu/tag/062Z}{Tag 062Z}]{stacks-project}).\footnote{Recall that a closed immersion $i \colon Z \to X$ of classical schemes which is locally of finite presentation (in the classical sense) is quasi-regular if and only if the classical conormal sheaf is a vector bundle which agrees with the classical conormal bundle \cite[\href{https://stacks.math.columbia.edu/tag/063M}{Tag 063M}]{stacks-project}. Theorem~\ref{ThmQ} and its proof can thus be seen as a strong realization of Kontsevich's  `hidden smoothness principle'.}

\textbf{The construction} of the \emph{derived extended Rees algebra} is in the present work induced from a pair of morphisms of derived algebraic contexts and is thus more robust than all previous constructions. We expect that our approach will be useful for constructing derived analogues of multi-centered Rees algebras. Such will likely  use the construction of the \emph{$\M$-filtered derived algebraic context} $\CCC^\M$ for a given ordered monoid $\M$ (Corollary~\ref{Cor:CM_DAC}). One can further iterate the procedure of replacing $\CCC$ by $\CCC^{\M}$ (possibly with different $\M$ at each step), and one question to ask is how this iterated procedure interacts with taking derived extended Rees algebras. 

\textbf{A potential application} of this story is towards a fully coherent Gysin map or purity transformation in motivic homotopy theory.
Namely, since the class of nonconnectively affine morphisms is now closed under the operation of taking the deformation space by Theorem~\ref{ThmF}, the current setup should also support a very general derived analogue of double (and higher) deformation spaces, as in \cite{AdeelFundamental}, and one should be able to describe these higher deformation spaces using derived multi-centered Rees algebras. 

\textbf{We close our notes} with a short generalization of the procedure of \emph{transmutation} from \cite{BhattF} in the setting of $\CCC_{\geq 0}$-stacks (\S \ref{Subsec:Trans}). The procedure starts with a $\CCC_{\geq 0}$-stack $\A$ valued in $\DAlg(\CCC)_{\geq 0}$, and produces for any stack $X$ a new prestack $X^{\A}$ defined by precomposing $X$ with $\A$. We generalize a result from \cite{Weil} and show how the deformation space of a given map $X \to Y$ of $\CCC_{\geq 0}$-stacks can be obtained through transmutation (Example~\ref{Exm:d2nb_via_transm}). 

One defines the \emph{$\A$-cohomology} of $X$ as the global sections of $X^{\A}$. We further include a first exploration of Rees algebra stacks induced from ideal stacks, and use this theory to give, for an ideal stack $\I \to \A$, a natural filtration 
    \[ H_{\B}(X) \to H_{\B_{(2)}}(X) \to \cdots \to H_{\A}(X) \]
on the $\A$-cohomology of $X$ by the $\B_{(k)}$-cohomology of $X$ for $\B_{(k)} \coloneqq \A/\I^k$ (Definition~\ref{Def:I_adic_filtr_coh}), which we show is exhaustive over $\CCC = \Mod_\Z$ when $\A$ is $\I$-adically complete (Proposition~\ref{Prop:exhaustive}).

\subsection*{Further connections to the literature}
The derived extended Rees algebra of $A \to B$ in $\DAlg(\Mod_\Z)$ coincides with the (unbounded) derived $I$-adic filtration of $B$ over $A$, where $I$ is the derived ideal associated to $A \to B$, up to identifying $\Z$-graded $A[t^{-1}]$-modules with $\Z$-filtered $A$-modules and their respective derived rings. For a definition of \emph{derived $I$-adic filtrations}, see \cite[Def.~2.5]{BrantnerFormal}. The comparison to derived extended Rees algebras was already observed in \cite[\S2.4]{MaoRevisiting}. The latter work primarily focuses on the $\CAlg$-valued version. It generalizes \cite[\S 4.2]{BhattCompletions} which is mostly over $\Q$, and it goes back to Quillen. The comparison in \cite[\S2.4]{MaoRevisiting} was observed with respect to the definition of the derived extended Rees algebra given in \cite{HekkingGraded}. Using our approach in \S \ref{subsec:Rees_Iadic} substantially simplifies the argument.

In \cite[Def.~9.1]{antieau2025filtrationscohomologyicrystallization} one finds the definition of the \emph{Hodge filtered derived infinitesimal cohomology}, which is the left adjoint of the functor
\[ \Fil^{\geq 0}(\Mod_A) \to \DAlg(\Mod_A) \colon F  \mapsto F^0/F^1, \]
where $A \in \DAlg(\Mod_\Z)$. Comparing universal properties, it is easy to see that this is also the derived Rees algebra. Our constructions thus give a globalization and generalization of this.

Let $A \in \DAlg(\Mod_\Z)$ be connective and let $I_0 \subset \pi_0(A)$ be finitely generated, say with generators $f_1,\dots,f_k$. Let $I$ be the derived ideal associated to the map $A \to A \sslash (f_1,\dots,f_k)$ (see Example \ref{Exm:classical_ideals}), and let $M \in \Mod_A$. 
Then $M$ is derived $I$-adically complete if and only if it is $I_0$-complete in the sense of \cite[Def.~4.2.1]{DAGXII}, by Proposition~\ref{Prop:derived_I_complete_T} and \cite[Prop.~4.2.7, Cor.~4.2.12]{DAGXII}. By Theorem~\ref{ThmB}, this confirms that in this case all notions of derived completeness coincide, c.f.\ \cite[\href{https://stacks.math.columbia.edu/tag/091S}{Tag 091S}]{stacks-project}.

This paper originated from a desire to use the tools developed herein to study the structure of prismatic cohomology. While we do not pursue them in this paper, such applications are a key focus of the closely related independent work \cite{SahaiDAPG}. Of relevance to this paper, using the conventions of the previous paragraph, \cite[Thm.~1.7]{SahaiDAPG} establishes the folklore result that there is a canonical equivalence of $\infty$-categories $\QCoh(\Spf(A))\simeq\Mod_A^{I\textup{-comp}}$ for $\Spf(A)$ the derived $I$-adic formal scheme associated to $A$ and $\Mod_A^{I\textup{-comp}}$ the localization of $\Mod_A$ spanned by derived $I$-complete $A$-modules. \cite[Thm.~1.9]{SahaiDAPG} then establishes a graded generalization of this result, followed by analysis of the globalized situation. This includes a generalization of the description of $\Theta$ as classifying stack for filtered objects.

We sketch two more intended applications of these notes, in which we fully adopt the language used in the main body of the text.

\begin{Fut**}[\textbf{Derived-geometric Nullstellensatz}]
    Let $X$ be a $\CCC_{\geq 0}$-stack. Write $\Ideal(X)_{\geq 0}$ for the category of connective, quasi-coherent ideals on $X$, and write $\Cld(X)$ for the category of closed immersions $Z \to X$. Taking vanishing loci furbishes an adjunction
    \[ (I \dashv V) \colon \Cld(X) \rightleftarrows \Ideal(X)_{\geq 0}^\op  , \]
    which in fact is an adjoint equivalence.\footnote{One could call this a `categorical Galois correspondence'.}

    Let $\RRR$ be the right Bousfield localization 
    \[ (i \dashv F) \colon \RRR \rightleftarrows \Cld(X) \]
    at all infinitesimal neighborhoods of the form $D \to D^{(n)}$, where $D \to X$ is a virtual Cartier divisor and $n \geq 1$. 

    Say that $Z \to X$ is \emph{reduced} when $iF(Z) \to Z$ is invertible. Say $\II \in \Ideal(X)_{\geq 0}$ is \emph{radical} if, whenever there is a map of the form $V(\II) \to D^{(n)}$ over $X$---where $D \to X$ is a virtual Cartier divisor and $n$ a natural number---then there is a map $V(\II) \to D$ over $X$ as well. The question is: for which $\CCC$ and which $X$ does it hold that every reduced $Z \in \Cld(X)$ is such that $I(Z)$ is radical, and for which $\CCC$ and which $X$ does the converse hold? Of course, part of the question is also to understand why this fails when it does.

    The answer should boil down to the classical Nullstellensatz when working over $\CCC = \Mod_k$ for an algebraically closed field $k \in \Mod_\Z$ and with affine $X$. 
\end{Fut**}

\begin{Fut**}[\textbf{Derived global analytic geometries}]
    \emph{Global analytic geometry} is a project that aims to give a unified approach to different flavors of analytic geometry in such a way that it also includes algebraic geometry. The term seems to originate from the eponymous work by Paugam, but the idea is much older with a history that includes Grothendieck, Berkovich, Zariski, Huber, and many others \cite{PaugamGlobalanalyticgeometry}. 
    Recent work in progress on analytic stacks by Clausen--Scholze can be understood as a continuation of this story.\footnote{This approach is based on their theory of \emph{condensed mathematics}, closely related to the independent theory of \emph{pyknotic objects} \cite{BarwickPyknotic}.}  At the time of this writing their work has not been formally written up yet, although the authors have made lecture series and notes widely available: traces can be found in \cite{KestingCategorical,AnschtzDescent,CamargoAnalytic,MannPadic}.
    Independently and culminating in a derived theory, there is the program developed by Bambozzi, Ben-Bassat, Kelly, Kremnizer, Mukherjee (in various configurations), see \cite{BenbassatPerspective,BenbassatAnaly,KellyAnalyticHKR,BenbassatNonarch,BambozziDagger}. This story revolves primarily around the derived category $\IndBan_R$ for a Banach ring $R$ and variations thereof such as bornological $R$-modules. A third strand can be found in the work of Holstein, Porta, Yue Yu (in various configurations), see \cite{ HolsteinAnalytification,PortaDerivedcomplexanalyticgeometryI, PortaDerived,PortaHigher,PortaDerivedcomplexanalyticgeometryII}. The latter focuses primarily on (specific flavors of) derived complex analytic and non-Archimedean analytic stacks, viewed through the lens of Lurie's (pre)geometries.

    The second-named author is currently developing a theory of derived global analytic geometries (joint with Oren Ben-Bassat and Jack Kelly). Roughly speaking, our approach has three steps. First, we give a slight generalization of the notion of derived algebraic contexts---called \emph{localized contexts}. Second, we define a notion of \emph{derived pre-analytic rings} over a fixed localized context $\EEE$, and make precise when such derived pre-analytic rings over $\EEE$ are connective and when they are normalized. We define an \emph{analytic $\EEE$-algebra} as a derived pre-analytic ring which is connective and normalized, and write $\AnAlg(\EEE)$ for the category of such. In the final step, we study stacks over $\AnAlg(\EEE)^\op$ which will be our \emph{derived analytic stacks over $\EEE$.} This approach can be seen as a synthesis of the condensed approach and the Banach-based approach to analytic geometry, together with a derived extension of the theory along the lines of \cite{RaksitDerived} for the affine picture and \cite{OrenJeroenblowups} for the global picture. The axiomatic framework morally provides a base which is deeper than the $\otimes$-unit in any given approach, and we expect to recover all known approaches to derived analytic geometry and many more by varying $\EEE$.

    A first application we have in mind with this is a unified treatment and vast generalization of HKR-type theorems. A long-term goal is to study the motivic homotopy theory of these new kinds of derived analytic stacks. We expect that the present work can be generalized to the setting of these derived analytic stacks over $\EEE$, and we expect that an appropriate version of adic completeness induces an example of a pre-analytic structure on $A$, for a given ideal pair $(A,I)$ over $\EEE$.
\end{Fut**}

\subsection*{Notation \& conventions}
From here on, everything will be implicitly $\infty$-categorical and derived, unless otherwise stated. We ignore `size issues' throughout. Our preferred solution is assuming an ample collection of universes in our meta-theory. 

For a category $\CCC$ and objects $x,y$ in $\CCC$ we write $\CCC(x,y)$ for the mapping space of morphisms $x \to y$ in $\CCC$.

For $M' \to M$ a map in a stable category, we may write the cofiber as $M/M'$.
For $K$ an algebra, module, space, spectrum, etc., we write $x \in K$ to mean $x \in \pi_0 K$.
For $\CCC$ a closed symmetric monoidal category, we write $\Map_\CCC(-,-)$ for the internal mapping object or simply $\Map(-,-)$ when $\CCC$ is clear from context.

Adjunctions will be written as
\[ (F \dashv G) \colon \CCC \rightleftarrows \DDD \]
where $F \colon \CCC \to \DDD$ is the left adjoint.\footnote{The notation is inspired by \cite{MacCategories}.}

When speaking of $\CCC_{\geq 0}$-stacks over a derived algebraic context $\CCC$, we say we are in the \emph{algebraic setting} to mean the case where $\CCC = \Mod_\Z$ and the Grothendieck topology on $\DAlg(\Mod_\Z)_{\geq 0}^\op$ is the \'{e}tale topology.

\subsection*{Acknowledgements}
The authors thank Zhouhang Mao for helpful and insightful discussions at multiple stages of writing these notes.
They also thank Edith H\"{u}bner for their interest in their work, and Universit\"{a}t M\"{u}nster for supporting a visit in August 2024.

The first named author was partially supported by the National Science Foundation (project number DMS-2200804). The second named author was supported by the Knut and Alice Wallenberg Foundation (project number 2021.0287), and by the DFG via the SFB 1085: Higher Invariants (project number 224262486). This work was also partially funded  by the Deutsche Forschungsgemeinschaft under Germany's Excellence Strategy EXC 2044 –390685587, Mathematics M\"{u}nster: Dynamics–Geometry–Structure.

\section{Background}
We give a reminder of the basic definitions surrounding derived rings due to Bhatt, Mathew, and Raksit in \S\ref{subsec:Derived_algebraic_contexts}, following \cite{RaksitDerived}. We also sketch the approach in \cite{OrenJeroenblowups} to stacks which are affine-locally the spectra of (connective) derived rings (\S\ref{subsec:geometric_cntxt}), including an overview of the deformation to the normal bundle and virtual Cartier divisors in this setting (\S\ref{subsec:Deformation_2NB} and \S\ref{subsec:VCD}). The latter follows and (partially) generalizes \cite{KhanVirtual, Weil}.
\subsection{Derived algebraic contexts}
\label{subsec:Derived_algebraic_contexts}
\begin{Def}
    A \emph{derived algebraic context} consists of
    \begin{itemize}
        \item a presentably symmetric monoidal, stable  category $\CCC$,
        \item endowed with a right complete $t$-structure
        $(\CCC_{\geq0},\CCC_{\leq0})$ compatible with the symmetric monoidal structure,
        \item and a symmetric monoidal full subcategory $\CCC^0\subseteq\CCC^{\heart}$, closed under finite coproducts in $\CCC$ and symmetric powers in $\CCC^{\heart}$, such that $\PPP_{\Sigma}(\CCC^0)\simeq\CCC_{\geq0}$.\footnote{Unpacking the requirement that $\CCC^0\subseteq\CCC^{\heart}$ is a a symmetric monoidal subcategory yields that $\pi_0(X_{\Sigma_n}^{\tensor n})\in\CCC^0$ for every $n\geq0$ and every $X\in\CCC^0$.}\footnote{Compatibility here means that $\CCC_{\leq0}$ is closed under filtered colimits, the unit object is connective, and the tensor product of connective objects is connective.}
    \end{itemize}
    A \emph{morphism of derived algebraic contexts} is a colimit-preserving symmetric monoidal functor $\CCC\to\DDD$ between derived algebraic contexts which is right $t$-exact and carries $\CCC^0$ into $\DDD^0$.\footnote{Recall that a functor is \emph{right $t$-exact} if it preserves connective objects.}
\end{Def}

From now on, let $\CCC$ denote a fixed derived algebraic context with unit $\mathrm{1}_\CCC$. Recall that $\LSym \colon \CCC \to \CCC$ is the monad obtained by taking the right-left extension of the symmetric algebra functor $\Sym_{\CCC^{\heartsuit}} \colon \CCC^0 \to \CCC$, and that $\DAlg(\CCC)$ is the category of $\LSym$-algebras in $\CCC$. We refer to objects of $\DAlg(\CCC)$ as \emph{derived $\CCC$-algebras} and objects of $\DAlg(\CCC)_{\geq0}:=\DAlg(\CCC)\times_{\CCC}\CCC_{\geq0}$ as \emph{connective $\CCC$-algebras}. This yields an adjunction
\[ (\iota \dashv \tau_{\geq 0}) \colon \DAlg(\CCC)_{\geq0} \rightleftarrows \DAlg(\CCC) \]
which commutes with the forgetful functors and the adjunction $\CCC_{\geq0} \rightleftarrows \CCC$ induced by the connective covers functor $\tau_{\geq0}(-)$ in the obvious way.

We write $\CAlg(\CCC)$ for the category of commutative algebra objects in $\CCC$. There is then a forgetful functor $\Theta \colon \DAlg(\CCC) \to \CAlg(\CCC)$, which associates to a derived $\CCC$-algebra $A$ the underlying commutative algebra $\Theta A$. Let $\Mod_A(\CCC)$ be the category of $\Theta A$-modules in $\CCC$, simply called \emph{$A$-modules}. Put $\DAlg_A(\CCC) \coloneqq \DAlg(\CCC)_{A/}$, and $\CAlg_A(\CCC) \coloneqq \CAlg(\CCC)_{\Theta A/}$. Endow $\DAlg_A(\CCC)$ with the coCartesian symmetric monoidal structure. Then the forgetful functors $\DAlg_A(\CCC) \to \CAlg_A(\CCC) \to \Mod_A(\CCC)$ are all symmetric monoidal.

A morphism of derived algebraic contexts $F \colon \CCC \to \DDD$ with right adjoint $G$ induces an adjunction on derived algebras, also written $(F \dashv G) \colon \DAlg(\CCC) \rightleftarrows \DAlg(\DDD)$. Both commute with the forgetful functors and $F$ with $\LSym$ in the obvious way.

\subsection{Graded objects}
Recall that $\CCC^\Z \coloneqq \Fun(\Z,\CCC)$ is the category of $\Z$-graded objects in $\CCC$, which is a derived algebraic context with pointwise $t$-structure. We let $M \mapsto M(0)$ be the fully faithful embedding $\CCC \to \CCC^\Z$ of algebraic contexts, which sends $M \in \CCC$ to the $\Z$-graded module which has $M$ concentrated in degree $0$.\footnote{By \emph{degree} we always mean \emph{homogeneous degree} with respect to a given grading (or, more generally, a given filtration).} 

For $N \in \CCC^\Z$ and $i \in \Z$ we write $N(i) \in \CCC^\Z$ for the \emph{twist of $N$ by $i$}, which is the graded object such that $N(i)^j \coloneqq N^{i+j}$ for all $j \in \Z$.\footnote{We are thus following the grading convention as found in \cite[\href{https://stacks.math.columbia.edu/tag/00JL}{Tag 00JL}]{stacks-project}.} Then $\CCC^\Z_{\geq 0}$ has as compact projective generators $N(i)$ for $N \in \CCC^0$ and $i \in \Z$. Observe that $N(i')(i) \simeq N(i'+i)$, and for $M \in \CCC$ we put $M(i) \coloneqq M(0)(i)$.

We write $\DAlg^\Z(\CCC) \coloneqq \DAlg(\CCC^\Z)$.

\subsection{Geometry over a given context}
\label{subsec:geometric_cntxt}
Write $\Aff(\CCC) \coloneqq \DAlg(\CCC)_{\geq 0}^{\op}$. In order to formulate the deformation to the normal bundle, we need to be able to do geometry relative to the context $\CCC$. To this end, endow $\Aff(\CCC)$ with a subcanonical topology $J$ such that the functor 
\[ \Mod_{(-)} \colon \DAlg(\CCC)_{\geq 0} \to \Cat \]
satisfies $J$-descent. This is part of the definition of a geometric context as defined in \cite{OrenJeroenblowups}. We review some key constructions from \textit{loc.\ cit.}

Let $\St(\CCC)$ be the  category of presheaves on $\Aff(\CCC)$ that satisfy $J$-descent. Objects of $\St(\CCC)$ are called \emph{$\CCC_{\geq 0}$-stacks} (or \emph{stacks} for short). Yoneda induces a fully faithful embedding
\[ \Spec(-) \colon \DAlg(\CCC)_{\geq 0}^{\op} \to \St(\CCC). \]
Restricted Yoneda also induces a functor
\[ \Spec^{\nc}(-) \colon \DAlg(\CCC)^{\op} \to \St(\CCC), \]
which in general is not fully faithful. 

For a stack $X \in \St(\CCC)$, the category of quasi-coherent modules 
\[ \QCoh(X) \simeq \lim_{\Spec(A) \to X} \Mod_A(\CCC) \]
is defined via right Kan extension of the functor $\Spec(A) \mapsto \Mod_A(\CCC)$ along $\Spec(-)$. Likewise the category of quasi-coherent algebras 
\[ \QAlg(X) \simeq \lim_{\Spec(A) \to X} \DAlg_A(\CCC)\]
is defined via  right Kan extension of $\Spec(A) \mapsto \DAlg_A(\CCC)$. Both of these categories canonically carry a symmetric monoidal structure, and $\QCoh(X)$ is stable and carries a $t$-structure. For a morphism $f \colon X \to Y$, we acquire an adjunction
\[ f^* \dashv f_* \colon \QCoh(Y) \rightleftarrows \QCoh(X), \]
and likewise on quasi-coherent algebras.

\begin{Not}
    If $\CCC$ is clear from context, we may write the categories $\DAlg_A(\CCC), \Aff(\CCC),\St(\CCC),...$ as $\DAlg_A, \Aff, \St,...$.
\end{Not}

Let $X$ be a stack. The functor $\LSym$ in the affine case induces an adjunction
\[ \LSym_{\OOO_X} \dashv U_X \colon \QCoh(X) \rightleftarrows \QAlg(X), \]
where $U_X$ is the forgetful functor. For $f \colon X \to Y$ and $M \in \QCoh(Y)$, it holds that $f^*\LSym_{\OOO_Y}(M) \simeq \LSym_{\OOO_Y}(f^*M)$. By passing to right adjoints, we see that the diagram
\begin{center}
    \begin{tikzcd}
        \QAlg(X) \arrow[r, "f_*"] \arrow[d, "U_X"] & \QAlg(Y) \arrow[d, "U_Y"] \\
        \QCoh(X) \arrow[r, "f_*"] &\QCoh(Y)
    \end{tikzcd}
\end{center}
commutes. 

We have the obvious $\Z$-graded analogues, which are written $\QCoh^\Z(-)$ and $\QAlg^\Z(-)$.

Stacks of the form $\Spec(A)$ for $A \in \DAlg(\CCC)_{\geq 0}$ are called \emph{affine}. 
A morphism $f \colon X \to Y$ of stacks is \emph{affine} if $X_A \coloneqq X \times_Y \Spec A$ is  affine, for all $\Spec A \to Y$.

\subsection{Deformation to the normal bundle}
\label{subsec:Deformation_2NB}
We recall the deformation to the normal bundle from \cite{OrenJeroenblowups} (which generalizes the relevant parts of \cite{KhanVirtual,Weil}). 

Let $\mathrm{1}_\CCC[t^{-1}]$ be the derived $\Z$-graded $\CCC$-algebra $\LSym(\mathrm{1}_\CCC(1)) \in \DAlg(\CCC^\Z)$, i.e., $\mathrm{1}_\CCC[t^{-1}]$ is the algebra freely generated by the variable $t^{-1}$ in degree $-1$. The unique morphism of derived algebraic contexts $\Mod_\Z \to \CCC$ extends uniquely to a morphism of derived algebraic contexts $\Mod^\Z_\Z \to \CCC^\Z$ which commutes with the various insertion and evaluation functors. Then $\mathrm{1}_\CCC[t^{-1}]$ is also the image of $\Z[t^{-1}]$ under the functor $\DAlg(\Mod^\Z_\Z) \to \DAlg(\CCC^\Z)$. 

 Likewise, we have the derived $\CCC$-algebra $\mathrm{1}_\CCC[t,t^{-1}] \in \DAlg(\CCC^\Z)$, which is the image of $\Z[t,t^{-1}]$ under the functor $\DAlg(\Mod^\Z_\Z) \to \DAlg(\CCC^\Z)$. Put $\A^1_\CCC \coloneqq \Spec(\mathrm{1}_\CCC[t^{-1}])$ and $\G_{m,\CCC} \coloneqq \Spec (\mathrm{1}_\CCC[t,t^{-1}])$. 

 \begin{Not}
     Let $\Theta_\CCC \coloneqq [\A_\CCC^1/\G_{m,\CCC}]$ be the quotient stack of the scaling action of $\G_{m,\CCC}$ on $\A^1_{\CCC}$ of degree $-1$. Write $i \colon B\G_{m,\CCC} \to [\A^1_{\CCC}/\G_{m,\CCC}]$ for the morphism of stacks induced by the zero section, and $j \colon [\G_{m,\CCC}/\G_{m,\CCC}] \simeq * \to [\A^1_\CCC/\G_{m,\CCC}]$ for the map induced by the canonical map $\G_{m,\CCC} \to \A^1_
    \CCC$.
 \end{Not}

As before, if $\CCC$ is clear from context, we may omit it from notation and write $\A^1$ for $\A^1_\CCC$ etc.

\begin{Def}
     The functor that sends a morphism $W \to X$ of $\CCC_{\geq 0}$-stacks to the \emph{deformation space} $\cD_{W/X} \in \St(\CCC)$ is right adjoint to the functor
    \begin{align*}
        \St_{\Theta} &\to \Arr(\St) \\
        T &\mapsto (T_0 \coloneqq T \times_{\Theta} B\G_{m} \to T)
    \end{align*}
    Put $D_{W/X} \coloneqq \cD_{W/X}\times_{B\G_{m}} *$, endowed with the canonical $\G_{m}$ action so that $[D_{W/X}/\G_{m}] \simeq \cD_{W/X}$. 
\end{Def}
By comparing universal properties, one obtains:
\begin{Lem}
    The functor $W \mapsto \cD_{W/X}$ is the Weil restriction of the morphism $W \times B\G_m \to X \times B\G_m$ along $X \times B\G_m \to X \times \Theta$.
\end{Lem}

Let $\D_1$ be the affine $\G_{m}$-scheme $\{0\} \times_{\A^1}  \{0\}$. 

\begin{Def}
    The \emph{normal bundle} $N_{W/X}$ of $(W \to X) \in \St_\CCC$ is the $\G_{m,\CCC}$-equivariant Weil restriction of $W \times \D_1 \to X \times \D_1$ along the natural map $z \colon X \times \D_1 \to X \times \{0\}$. Put $\cN_{W/X} \coloneqq [N_{W/X}/\G_m]$.
\end{Def}

Write $s \colon X \to X \times \D_1$ for the diagonal, and let $\epsilon \colon z^*z_*(W \times \D_1) \to W \times \D_1$ be the counit. Then pulling back $\epsilon$ along $s$ induces a morphism
\[ s^*(\epsilon) \colon N_{W/X} \to W. \]

\begin{Prop}
\label{Prop:D2NB}
We have a $\G_{m,\CCC}$-equivariant Cartesian diagram
		\begin{center}
		\begin{tikzcd}
			W\ar[r]\ar[d] & W\times \A^1_{\CCC}\ar[d] & W\times \G_{m}\ar[l]\ar[d]\\
			N_{W/X}\ar[r]\ar[d] & D_{W/X}\ar[d] & X\times \G_{m}\ar[l]\ar[d]\\
			X\ar[r] & X\times \A^1_{\CCC} & X\times \G_{m}\ar[l]		
		\end{tikzcd}	
	\end{center}
        natural in $(W \to X) \in \St_\CCC$. If $W \to X$ admits a cotangent complex $L_{W/X}$, then 
	\[ N_{W/X} \simeq \widetilde{\Spec} (\LSym(L_{W/X}[-1])).  \]
\end{Prop}

\begin{proof}
    The middle column of the Cartesian diagram is induced by naturality of $D_{(-)/(-)}$. Base-change for Weil restrictions then induces the desired squares and shows that they are Cartesian. The last claim is \cite[Thm.~4.15]{OrenJeroenblowups}.
\end{proof}

\subsection{Virtual Cartier divisors}
\label{subsec:VCD}
Let $G \in  \St(\CCC)$ be an affine group stack. Recall that $\QCoh(BG)$ is equivalent to the category of comodules over $\OOO_G$ in $\CCC$. It follows that $\QCoh(B\G_m) \simeq \CCC^\Z$ as symmetric monoidal categories \cite[Prop.~2.25]{OrenJeroenblowups}. With a similar argument, it follows that $\QCoh(\Theta) \simeq \Mod_{\mathrm{1}_\CCC[t^{-1}]}(\CCC^\Z)$ as symmetric monoidal categories. One can show that these equivalences are $t$-exact.

\begin{Def}
    A \emph{virtual Cartier divisor} is a morphism $D \to T$ of stacks such that there is a Cartesian diagram
    \begin{center}
        \begin{tikzcd}
            D \arrow[r] \arrow[d] & T \arrow[d] \\
            B\G_m \arrow[r] & \Theta.
        \end{tikzcd}
    \end{center}
    Let $W \to X$ be a morphism of stacks. Then a \emph{virtual Cartier divisor over $W \to X$} is a commutative diagram
    \begin{center}
        \begin{tikzcd}
            D \arrow[r] \arrow[d] & T \arrow[d] \\
            W \arrow[r] & X
        \end{tikzcd}
    \end{center}
    in $\St(\CCC)$ such that $D \to T$ is a virtual Cartier divisor.
\end{Def}

Let $W \to X$ be a morphism of stacks. Observe that $\cN_{W/X} \to \cD_{W/X}$ naturally lives over $W \to X$.

\begin{Prop}
    The morphism $\cN_{W/X} \to \cD_{W/X}$ is the universal virtual Cartier divisor over $W \to X$, in the sense that any virtual Cartier divisor $D \to T$ over $W \to X$ is obtained as a pullback from the universal one. 
\end{Prop}

\begin{proof}
    This follows from the deformation to the normal bundle.
\end{proof}

\section{Generalized filtered objects} 
We provide a general machinery for producing new derived algebraic contexts from the given derived algebraic context $\CCC$ and an ordered (classical) monoid $\M$. Our Corollary~\ref{Cor:CM_DAC} generalizes \cite[Lem.~2.18]{OrenJeroenblowups} where $\M$ is required to have no nontrivial morphisms, and it unifies the latter approach with the construction of the derived algebraic context of $\M'$-filtered objects in $\CCC$ for $\M' = \{0<1\}, \N,\Z$. The case $\M' = \{0<1\}$ underpins the theory of ideal pairs over $\CCC$ as explained in \S\ref{Sec:Ideal_affine}.
\subsection{First definitions}
\label{subsec:first_defs}
An \emph{ordered monoid} is a symmetric monoidal, thin category $\M$. For $i,j \in \M$ we write $i \leq j$ if there is an arrow $j \to i$. The monoidal structure on $\M$ is written additively. 
For a presentably symmetric monoidal, stable category $\DDD$, write
\[ \DDD^\M \coloneqq \Fun(\M,\DDD), \]
endowed with the Day convolution symmetric monoidal structure. Recall that $\DDD^\M$ is again presentably symmetric monoidal and stable. 
An object $X \in \DDD^\M$ is called an \emph{$\M$-filtered} object in $\DDD$, and we write $i \mapsto X^i$ for $i \in \M$. 
A \emph{morphism of ordered monoids} is a lax monoidal functor $\alpha \colon \M \to \K$ between ordered monoids. Given such, we obtain adjunctions
\[ (\alpha_! \dashv \alpha^*) \colon \DDD^\M \rightleftarrows \DDD^\K \]
where $\alpha^* \colon \DDD^\K \to \DDD^\M$ is restriction along $\alpha$ and $\alpha_!$ is left Kan extension. 

\begin{Prop}
\label{Prop:alpha_monoidal}
    The functor $\alpha^* \colon \DDD^\K \to \DDD^\M$ canonically carries a lax monoidal structure, hence $\alpha_!$ carries an oplax monoidal structure. Moreover, if $\alpha$ is strong then so is $\alpha_!$.
\end{Prop}

\begin{proof}
    This follows from the universal property of the Day convolution symmetric monoidal structure. An explicit reference is \cite[Cor.~3.8]{NikolausYoneda}. 
    
    Informally, we can also see this from the pointwise formula for the Day convolution symmetric monoidal structure. Namely, since $\alpha$ is lax symmetric monoidal, there is a natural transformation
    \[ \tau \colon (\mu_\K \circ (\alpha,\alpha) \to \alpha \circ \mu_{\M} ) \colon \M \times \M \to \K, \]
    where $\mu_\M \colon \M \times \M \to \M$ and $\mu_\K \colon \K \times \K \to \K$ are the monoidal opertions on $\M$ and on $\K$ respectively. Hence we have the Beck--Chevalley map
    \[ {\mu_\M}_! (\alpha,\alpha)^* \xrightarrow{\eta} {\mu_\M}_! (\alpha,\alpha)^*\mu_\K^*{\mu_\K}_! \xrightarrow{\tau} {\mu_\M}_! \mu_\M^* \alpha^* {\mu_\K}_! \xrightarrow{\epsilon} \alpha^* {\mu_\K}_!, \]
    which exhbits a lax monoidal structure on $\alpha^*$ (omitting the higher coherency data). The second claim follows by naturality of left Kan extensions.
\end{proof}

\begin{Rem}
\label{Rem:ins}
    For $i \in \M$ we have an adjunction 
    \[ (\ins^i \dashv \ev^i) \colon \DDD \rightleftarrows \DDD^\M, \]
    where $\ev^i \colon \DDD^\M \to \DDD$ is the functor which sends $Y \in \DDD^\M$ to the \emph{evaluation} $Y^i$ of $Y$ at $i$. The left adjoint---called the \emph{insertion} functor---is fully faithful and sends $X \in \DDD$ to the filtered object given explicitly by the formula
    \begin{equation*}
        \ins^i(X)^j:=
        \begin{cases}
            0, & j > i, \\
            X, & i\geq j,
        \end{cases}
    \end{equation*}
    for $j \in \M$, with obvious structure maps.

    For $\alpha \colon \M \to \K$ a morphism of ordered monoids it holds $\ev^i\alpha^* \simeq \ev^{\alpha(i)}$, hence $\alpha_!\ins^i \simeq \ins^{\alpha(i)}$.
\end{Rem}

\begin{Exm}
\label{Ex:ins_tensor}
    Consider the monoidal composition $\mu \colon \M \times \M \to \M$ as morphism of ordered monoids, and let $s,t\in \M$ and $X,Y \in \DDD$. By Remark~\ref{Rem:ins} it holds that
    \[ \ins^{s+t}(X \otimes Y) \simeq \mu_!\ins^{(s,t)}(X \otimes Y) \simeq \ins^s(X) \otimes \ins^t(Y), \]
    where the second equivalence is exactly the description of the Day convolution symmetric monoidal structure via left Kan extension along $\mu$. Consequently, we have $\ins^s(X)^{\otimes n}_{\Sigma_n} \simeq \ins^{n  s}(X^{\otimes n}_{\Sigma_n})$, where we have used that $\ins^{ns}(-)$ is a left adjoint, hence commutes with taking orbits $(-)_{\Sigma_n}$.
\end{Exm}

\begin{Lem}
    Suppose that $\DDD$ is endowed with a t-structure. Then there is a t-structure on $\DDD^\M$---called the \emph{neutral} t-structure---such that $X \in \DDD^\M$ is (co)connective if and only if $X^i$ is (co)connective for all $i \in \M$. Moreover:
    \begin{enumerate}
        \item If the $t$-structure on $\DDD$ is right-complete, then so is the neutral t-structure on $\DDD^\M$.
        \item If the $t$-structure on $\DDD$ is compatible, then so is the neutral t-structure on $\DDD^\M$.
    \end{enumerate}
\end{Lem}
\begin{proof}
    This is straightforward.
\end{proof}

\subsection{$\M$-filtered derived algebraic contexts}
Let still $\M$ be an ordered monoid, and consider the derived algebraic context $\CCC$. Write $\CCC^{\M,0}$ for the full subcategory of $\CCC^{\M,\heartsuit}$ spanned by finite coproducts of objects of the form $\ins^i(X)$ for $i \in \M$ and $X \in \CCC^0$.
\begin{Cor}
\label{Cor:CM_DAC}
    The category $\CCC^\M$ endowed with the neutral $t$-structure and the Day convolution symmetric monoidal structure is a derived algebraic context with $\CCC^{\M,0}$ as compact projective generators of $\CCC^\M_{\geq 0}$. 
\end{Cor}
\begin{proof}
     By Example~\ref{Ex:ins_tensor} and since the insertion functors are $t$-exact, the category $\CCC^{\M,0}$ is a symmetric monoidal subcategory of $\CCC^{\M,\heartsuit}$ which is closed under $\CCC^{\M,\heartsuit}$-symmetric powers. Clearly, $\CCC^{\M,0}$ is closed under finite coproducts. Hence it suffices to show that the left derived functor
    \[ F \colon \PPP_\Sigma(\CCC^{\M,0}) \to \CCC^\M_{\geq 0} \]
    of the inclusion $f \colon \CCC^{\M,0} \to \CCC^\M_{\geq 0}$ is an equivalence. Since the evaluation functors preserve sifted colimits, the insertion functors preserve compact projective objects, hence the fully faithful functor $f$ lands in compact projective objects. It remains to show that the image of $F$ generates $\CCC^\M_{\geq 0}$ under colimits \cite[Prop.~5.5.8.22]{LurieHTT}.

    Consider $\CCC^\M_{\geq 0}$ as full subcategory 
    \[ \CCC^\M_{\geq 0} \triangleq \Fun(\M,\CCC_{\geq 0}) \subset \PPP(\M^\op \times \CCC^0) \]
    via currying and the fact that $\CCC_{\geq 0} \simeq \PPP_\Sigma(\CCC^0) \subset \PPP(\CCC^0)$. For $i \in \M$ and $X \in \CCC^0$ it holds that $\ins^i(X)$ corresponds to the presheaf on $\M^\op \times \CCC^0$ representable by $(i,X)$. It follows that any $Y \in \CCC^\M_{\geq 0}$ is a colimit of objects of the form $\ins^i(X)$ with $i\in \M$ and $X \in \CCC^0$, from which the claim follows.
\end{proof}
\begin{Cor}
\label{Cor:shriek_cntxt}
    Let $\alpha \colon \M \to \K$ be a lax monoidal functor of ordered monoids and $\CCC$ a derived algebraic context. If $\alpha$ is strong, then $\alpha_! \colon \CCC^\M \to \CCC^\K$ is a morphism of derived algebraic contexts.
\end{Cor}
\begin{proof}
    Since $\alpha^*$ is $t$-exact the left adjoint $\alpha_!$ is right $t$-exact. Since $\alpha_!\ins^i(X) \simeq \ins^{\alpha(i)}(X)$ by Example~\ref{Ex:ins_tensor} it follows that $\alpha_!(\CCC^{\M,0}) \subset \CCC^{\K,0}$. Clearly $\alpha_!$ preserves colimits, and it is symmetric monoidal by Proposition~\ref{Prop:alpha_monoidal}.
\end{proof}
\begin{Exm}
    For $\alpha \colon \{0\} \to \M$ the canonical inclusion it holds that $\alpha_! \dashv \alpha^*$ is the adjunction $\ins^0 \dashv \res^0$. In particular, $\ins^0$ is a morphism of derived algebraic contexts. 
\end{Exm}

\subsection{The $\N$- and $\Z$-filtered case}
To emphasize the structure of an ordered monoid we write it as $(\M,\leq,+)^\op$ with the convention that there is an arrow $j \to i$ if and only if $i \leq j$. 
\begin{Def}
    For a derived algebraic context $\CCC$, write
     \begin{align*}
         \Fil(\CCC) \coloneqq \CCC^{(\Z,\leq,+)^\op} &&
          \Fil^{\geq 0}(\CCC) \coloneqq \CCC^{(\N,\leq,+)^\op} 
     \end{align*}
    for the derived algebraic contexts as constructed via Corollay~\ref{Cor:CM_DAC}.
\end{Def}
Precomposing with the functor $\Z \to (\Z, \leq)^\op$ yields a functor $V \colon \Fil(\CCC) \to \CCC^\Z$, which sends $X \in \Fil(\CCC)$ to the $\Z$-graded module $\bigoplus_i \ev^i(X)$. This has a left adjoint
\[ F \colon \CCC^\Z \to \Fil(\CCC) \]
which is a morphism of derived algebraic contexts by Corollary~\ref{Cor:shriek_cntxt}. 
Explicitly, $F$ sends $N \in \CCC^\Z$ to the filtered object $FN$ for which $\ev^i(IN) \simeq \bigoplus_{j \geq i} N^j$. 

We also have the functor
\[ \gr \colon \Fil(\CCC) \to \CCC^\Z \]
such that $\gr(X)^i \coloneqq \cofib(X^{i+1} \to X^i)$ for $X \in \Fil(\CCC)$. This functor is symmetric monoidal, and the composition
\[ \CCC^\Z \xto{F} \Fil(\CCC) \xrightarrow{\gr} \CCC^\Z \]
is equivalent to the identity \cite[Prop.~2.39]{MaoRevisiting}. 

\section{Ideal pairs: affine case}
\label{Sec:Ideal_affine}
We continue with the derived algebraic context $\CCC$. We define  ideal pairs over $\CCC$, and construct a correspondence between ideals $I \to A$ and morphisms of derived $\CCC$-algebras $A \to B$  in \S\ref{Subsec:Derived_ideal_pairs}. For an ideal $I \to A$ we show how $I$-adic filtrations can be constructed from the extended Rees algebra in \S\ref{subsec:Rees_Iadic}. In \S\ref{subsec:Pair_base} we relativize this story over a base $A \in \DAlg(\CCC)$, which will serve as a basis for the global counterpart in the next section.

\subsection{Derived ideal pairs}
\label{Subsec:Derived_ideal_pairs}
Consider the ordered monoid $([0,1],\leq,\max)^\op$.
\begin{Def}
    Write
    \[ \Pair(\CCC) \coloneqq \CCC^{([0,1],\leq,\max)^\op} \]
    for the derived algebraic contexts as constructed via Corollay~\ref{Cor:CM_DAC}.
\end{Def}
Explicitly, for $F,G\in\Pair(\CCC)$, the tensor product $F\tensor G$ is given by
$$(F^0\tensor G^1)\sqcup_{F^1\tensor G^1}(F^1\tensor G^0)\to F^0\tensor G^0$$
arising from the commutative square
\begin{center}
    \begin{tikzcd}
        F^1\tensor G^1 \arrow[r] \arrow[d] & F^0\tensor G^1 \arrow[d] \\
        F^1\tensor G^0 \arrow[r] & F^0\tensor G^0
    \end{tikzcd}
\end{center}
We call this the \emph{pushout product} symmetric monoidal structure. 

Observe that there is an equivalence of categories
\begin{align*}
    C \colon \Fun(([0,1],\leq)^\op,\CCC) &\to \Fun(([-1,0],\leq)^\op,\CCC) \\
    (X^1 \to X^0) &\mapsto (X^0 \to X^0/X^1)
\end{align*}
which sends an arrow $X^1 \to X^0$ to the cofiber $X^0 \to X^0/X^1$. The inverse is induced by taking fibers. 
\begin{Def}
\label{Def:ArC}
    Write $\Ar(\CCC)$ for the derived algebraic context with underlying category $\Fun(([-1,0],\leq)^\op,\CCC)$, induced by transport of structure via the equivalence $C$. We call the induced $t$-structure on $\Ar(\CCC)$ the \emph{ideal $t$-structure}.
\end{Def}
Clearly, the ideal $t$-structure on $\Ar(\CCC)$ is such that $f \colon Y^0 \to Y^{-1}$ is connective if and only if $Y^0,Y^{-1} \in \CCC_{\geq 0}$ and $f$ is surjective on $\pi_0$.
\begin{Prop}
    The induced symmetric monoidal structure on $\Ar(\CCC)$ is the pointwise symmetric monoidal structure.
\end{Prop}
\begin{proof}
    This is \cite[Prop.~3.1]{MaoRevisiting}.
\end{proof}

\begin{Def}
    A \emph{(derived) ideal pair (in $\CCC$)} is an object of $\DAlg(\Pair(\CCC))$. We will typically denote such an object $I\to A$ by $(A,I)$ and say that $I$ is an \emph{ideal} of $A$. We let $\DPair(\CCC)$ denote the category of derived ideal pairs in $\CCC$ and $\DPair_{\geq0}(\CCC) \coloneqq \DAlg(\Pair(\CCC))_{\geq 0}$ the full subcategory of connective ideal pairs. 
\end{Def}

Write $\Arr(\DDD) \coloneqq \Fun(\Delta^1,\DDD)$ for the arrow category in any category $\DDD$. Observe that there is a canonical equivalence
\[ \DAlg(\Ar(\CCC)) \simeq \Arr(\DAlg(\CCC)). \]

\begin{Not}
\label{Not:AmodI}
    For $(A,I) \in \DPair(\CCC)$, we write the corresponding object in $\DAlg(\Ar(\CCC))$ as $A \to A/I$. Conversely, for $(A \to B) \in \DAlg(\Ar(\CCC))$, we write the corresponding ideal pair as $(A, I(B))$.
\end{Not}

\begin{Exm}
\label{Exm:classical_ideals}
    Consider the context $\CCC = \Mod_\Z$, let $A \in \DAlg(\Mod_\Z)_{\geq 0}$ and $f_1,\dots,f_n \in A$. We define
    \[ A \sslash (f_1,\dots,f_n) \coloneqq A \otimes_{\Z[f_1,\dots,f_n]} \Z \]
    where $\Z[f_1,\dots,f_n] \to \Z$ sends $f_i$ to $0$.

    Now suppose $A$ is discrete. Then any ideal $I \subset A$ in the classical sense gives rise to an ideal pair $(A,I)$ in a canonical way. The associated quotient $A \to A/I$ is the classical quotient, as can be see from considering the associated long exact sequence.

    In particular we have the classical ideal $(f_1,\dots,f_n) \subset A$, and the  corresponding $A$-algebra $A/(f_1,\dots,f_n)$ generally differs from $A \sslash (f_1,\dots,f_n)$. Indeed, for a single element $f \in A$ we have an ideal pair $(A,fA)$ corresponding to the classical ideal $(f) \subset A$, and second we have the ideal pair $(A,I(f))$ corresponding to $A \to A\sslash(f)$.  Then $I(f) = A$ and the structure map $I(f) \to A$ is the multiplication map $\times f \colon A \to A$. There is a canonical map $(A,I(f)) \to (A,fA)$ and the induced map $A/I(f) \to A/fA$ on $A$-algebras identifies $A/fA$ with $\pi_0(A/I(f))$. The fiber of $A/I(f) \to A/fA$ is the cofiber of $\times f \colon A \to fA$, which is $\ker(f)[1]$. Thus, the difference between these two ideal pairs is captured precisely by the $f$-torsion in $A$.
\end{Exm}

\begin{Exm}
    By Corollary~\ref{Cor:CM_DAC} the compact projective generators of $\Pair(\CCC)_{\geq 0}$ are of the form $A \to A \oplus B$ where $A,B \in \CCC^0$ and the map is the natural inclusion. Consequently, $\Ar(\CCC)_{\geq 0}$ has compact projective generators of the form $A \oplus B \to B$, which give compact projective generators for $\DAlg(\Ar(\CCC))_{\geq 0}$ upon applying $\LSym_{\CCC}$ pointwise. For $\CCC = \Mod_\Z$ this recovers \cite[Thm.~3.23]{MaoRevisiting}.
\end{Exm}

\subsection{Rees algebras and $I$-adic filtrations}
\label{subsec:Rees_Iadic}
Consider the natural restriction functors
\[ \Fil(\CCC)\xrightarrow{(-)^{\geq 0}}\Fil^{\geq0}(\CCC)\xrightarrow{\res}\Pair(\CCC), \]
obtained by restriction along the morphisms of ordered monoids
\[ ([0,1],\leq,\max)^\op \to (\N,\leq, +)^\op \to (\Z,\leq, +)^\op. \]
Write $(-)^\ext \colon \Fil^{\geq 0}(\CCC) \to \Fil(\CCC)$ for the left adjoint of $(-)^{\geq 0}$.

\begin{Prop}
\label{Prop:ext_functors}
    The extension functor $(-)^\ext$ and the restriction functor $\res$ are morphisms of derived algebraic contexts.
\end{Prop}

\begin{proof}
    For $(-)^\ext$ this follows from Corollary~\ref{Cor:shriek_cntxt}. 
    
    For the second claim, it suffices to show that $\res$ is symmetric monoidal. This is straightforward, using that for filtered objects $X,Y$ it holds
    \[ (X \otimes Y)(n) \simeq \colim_{i + j \geq n} X(i) \otimes Y(j),\]  together with a cofinality argument. 
\end{proof}

\begin{Rem}
\label{Rem:psi_varphi_tilde}
    Let $\psi \colon \DDD \to \DDD'$ be a morphism of derived algebraic contexts which preserves limits and thus admits a left adjoint $\varphi$. Then, since the forgetful functors $\DAlg(\DDD)\to\DDD$ and $\DAlg(\DDD') \to \DDD'$ preserve limits and are conservative, we have that $\psi: \DAlg(\DDD) \to \DAlg(\DDD')$ preserves limits and so admits a left adjoint $\tilde{\varphi}$.
\end{Rem}

Since limits in functor categories are computed pointwise, combining Proposition~\ref{Prop:ext_functors} with Remark~\ref{Rem:psi_varphi_tilde} yields adjunctions
\begin{center}
    \begin{tikzcd}
        \DAlg(\Fil(\CCC)) \arrow[r, bend right, "(-)^{\geq 0}"{name=g}] & \DAlg(\Fil^{\geq 0}(\CCC)) \arrow[l, bend right, "(-)^{\ext}"{name=e}, swap] \arrow[phantom, from=e, to=g, "\dashv"  rotate=-90] \arrow[r, bend right, "\res"{name=s}]  & \DPair(\CCC). \arrow[l, bend right, "R"{name=R}, swap]
        \arrow[phantom, from=R, to=s, "\dashv" rotate=-90]
    \end{tikzcd}
\end{center}

\begin{Def}
\label{Def:Rees}
    The \emph{extended Rees algebra} of $(A\to B) \in \DAlg(\Ar(\CCC))$ is the $\Z$-graded algebra
    \[R^{\ext}_{B/A} \coloneqq R(A,I(B))^{\ext} \in \DAlg(\Fil(\CCC)) \simeq \DAlg(\Mod^\Z_{\mathrm{1}_\CCC[t^{-1}]}).\]
    We consider this an $A[t^{-1}]$ algebra via the natural map
    \[ A[t^{-1}] \simeq R^{\ext}_{A/A} \to R^{\ext}_{B/A} \]
    The \emph{Rees algebra} is the graded $A$-algebra $R_{B/A} \coloneqq (R^\ext_{B/A})^{\geq 0}$.
\end{Def}

We will compare the definition above with \cite{Weil} in \S \ref{Subsec:DeftoNB}.

\begin{Rem}
    The filtered algebra associated to the Rees algebra $R_{B/A}$ is denoted $p_!(A,I(B)) = I(B)^{\star}A$ in \cite[Not.\ 2.4.9]{HolemanDerivedDelta} and called the \emph{$I$-adic filtration} on $A$.
\end{Rem}

Let $\coins^{\geq 0} \colon \Pair(\CCC) \to \Fil^{\geq 0}(\CCC)$ be the right adjoint of $\res$, which exists since $\res$ is a morphism of derived algebraic contexts. Observe that $\coins^{\geq 0}(M^1 \to M^0)$ is the filtered object
\[ \cdots \xrightarrow{\id} M^1 \xrightarrow{\id} M^1 \to M^0. \]

\begin{Lem}
\label{Lem:Rees_fully_faithful}
    The functors $R$ and $(-)^{\ext}$ are fully faithful.
\end{Lem}

\begin{proof}
For fully faithfulness of $R$, it suffices to show that the induced functor $\coins^{\geq 0} \colon \DPair(\CCC) \to \DAlg(\Fil^{\geq 0}(\CCC))$ is fully faithful, which is clear. 

For fully faithfulness of $(-)^\ext$ it suffices to check this on the level of filtered modules. The latter follows from the fact that the morphism of ordered monoids $(\N,\leq,+)^\op \to (\Z,\leq, +)^\op$ is fully faithful
\end{proof}

Given a ideal pair $(A,I)$ and $n\geq 0$, we want to associate a new ideal pair $(A,I^n)$. To this end, let $q_n \colon \N \to \N$ be the dilation map $i \mapsto in$. Since this is a strong morphism of ordered monoids, we obtain a morphism of derived algebraic contexts
\[ {q_n}_! \colon \Fil^{\geq 0}(\CCC) \to \Fil^{\geq 0}(\CCC) \]
with right adjoint ${q_n}_! \dashv q_n^*$ such that $\ev^i q_n^* \simeq \ev^{in}$. As always, we also write ${q_n}_! \dashv q_n^*$ for the induced adjunction on derived $\Fil^{\geq 0}(\CCC)$-algebras.

\begin{Not}
    In what follows, $(A,I)$ is an ideal pair.
\end{Not}

\begin{Def}
    For $n\geq 0$ define the ideal pair $(A,I^n) \coloneqq \res q_n^*R(A,I)$.
\end{Def}

\begin{Def}
\label{Def:I_adic_filtration}
    Given $M\in\Mod_A$ and $n\geq0$, we define 
    \begin{align*}
        R(M,I) &\coloneqq R(A,I)\tensor_{R(A,0)} \ins^0(M), \\
        I^nM &\coloneqq \ev^n(R(M,I)).
    \end{align*}
    We call $R(M,I)$ the \emph{$I$-adic filtration on $M$}. These constructions define functors 
    \begin{align*}
        R(-,I) &\colon \Mod_A \to \Mod_{R(A,I)}(\Fil^{\geq0}(\CCC)) \\
        I^n &\colon \Mod_A\to\Mod_A.
    \end{align*}
    For $M=A$ we simply put $I^n \coloneqq I^nA$.
\end{Def}

\begin{Rem}
    There are now two meanings of $I^n$. First, via the ideal pair $(A,I)$, we have
    \[ I^n \triangleq \ev^n(R(A,I)\tensor_{R(A,0)} \ins^0(A)) \simeq \ev^nR(A,I) \simeq \ev^1 q_n^*R(A,I) , \]
    while via the ideal pair $(A,I^n)$ we have
    \[ I^n \triangleq  \ev^1(R(A,I^n)\tensor_{R(A,0)} \ins^0(A)) \simeq \ev^1R(A,I^n). \]
    Luckily these agree since the canonical map
    \[ \epsilon \colon R(A,I^n) \to q_n^* R(A,I) \]
    induced from the adjunction $R \dashv \res$ is an equivalence after applying $\res$, since $R$ is fully faithful by Lemma~\ref{Lem:Rees_fully_faithful}.
\end{Rem}

\begin{Rem}
    From the formula of the Day convolution symmetric monoidal structure via left Kan extension it follows that $I^nM$ is simply $I^n \otimes_A M$.
\end{Rem}

\begin{Warn}
    In general it does not hold $I^nI \simeq I^{n+1}$, where the second $I$ is the underlying $A$-module. Indeed, put $B \coloneqq A/I$, let $N_{B/A} \coloneqq L_{B/A}[-1]$ be the shifted cotangent complex, and let $F$ be the fiber of the natural map $B \otimes_A I \to N_{B/A}$. Assume $A,I,B$ are connective. As we shall later see, the deformation to the normal bundle implies $N_{B/A} \simeq I/I^2$ (Corollary~\ref{Cor:NB_is_NC}). Using the fiber sequence $I \otimes_A I \to I \to B \otimes_A I$, we obtain a Cartesian diagram
    \begin{center}
        \begin{tikzcd}
            I \otimes_A I \arrow[d] \arrow[r] & 0 \arrow[d] \\
            I^2\arrow[r] \arrow[d] & F \arrow[d] \arrow[r] & 0 \arrow[d]\\
            I \arrow[r] & B \otimes_A I \arrow[r] & I/I^2,
        \end{tikzcd}
    \end{center}
    and thus a fiber sequence
    \[ I \otimes_A I \to I^2 \to F. \]
    In particular, the difference between $I \otimes_A I \to I^2$ is measured by the Hurewicz map $B \otimes_A I[1] \to L_{B/A}$ over $\CCC = \Mod_\Z$, which in general is not an equivalence \cite[\S 25.3.6]{LurieSpectral}. For example, take $I,A,B$ discrete, and use that $\pi_1\LL_{B/A} \simeq I/I^2$  \cite[\href{https://stacks.math.columbia.edu/tag/08RA}{Tag 08RA}]{stacks-project}.
\end{Warn}

\begin{Warn}
\label{Warn:dilation_isnt_Rees}
    In general the canonical map 
    \[ \epsilon \colon R(A,I^n) \to q_n^* R(A,I) \]
    is not invertible. For a concrete example, consider the derived algebraic context $\CCC = \Mod_\Z$. Put $A \coloneqq \C[x,y]$ and let $I$ be the classical ideal $(x,y) \subset A$ so that $A/I \simeq \C$. Then $A \to A/I$ is quasi-smooth hence $R(A,I)$ is the classical Rees algebra of $A \to A/I$ so that the ideal powers $I^n$ are the classical powers of $I$. In particular, $I^2$ is the classical ideal $(x^2,xy,y^2)$. Put $B \coloneqq A/I^2$. Then $A \to B$ is no longer quasi-smooth since 
    \[ B \simeq  A / (x^2,xy,y^2) \]
    is the classical quotient by Example~\ref{Exm:classical_ideals}.

    Now suppose that $\epsilon \colon R(A,I^2) \to q_2^*R(A,I)$ is invertible. By passing to the associated graded and using Corollary~\ref{Cor:NB_is_NC}, we obtain an equivalence
    \[ L_{B/A}[-1] \simeq I^2/I^4 \]
    of $A$-modules, where $I^2 / I^4$ is the cofiber of the inclusion $I^4 \to I^2$ on the ideal powers which are classical. In particular, we would have that $L_{B/A}[-1]$ is discrete, and an explicit computation shows that this is not the case.\footnote{This follows from the fact that there are non-trivial $A$-linear relations between the generators $x^2,xy,y^2$ of $I^2$, which implies that $\pi_2\LL_{B/A} \not=0$  \cite[\href{https://stacks.math.columbia.edu/tag/09AM}{Tag 09AM}]{stacks-project}.} In fact, the natural $B$-module structure on $L_{B/A}[-1]$ does not carry over to $I^2/I^4$.
\end{Warn}

\subsection{Ideal pairs relative to a base}
\label{subsec:Pair_base}
Let $A \in \DAlg(\CCC)_{\geq 0}$ be given. Note that canonically $(A,0) \in \DAlg(\Pair(\CCC))_{\geq 0}$. As is the case for any derived algebraic context, the category
\[ \Pair_A(\CCC) \coloneqq \Mod_{(A,0)}(\Pair(\CCC)) \]
carries a canonical symmetric monoidal structure, and the $\LSym_{\Pair(\CCC)}$-monad induces a monadic adjunction
\[ (\LSym_{(A,0)} \dashv U_{(A,0)}) \colon \Pair_A(\CCC) \rightleftarrows \DPair_A(\CCC) \]
where $\DPair_A(\CCC) \coloneqq \DPair(\CCC)_{(A,0)/}$ and $U_{(A,0)}$ is the forgetful functor. Moreover, the category $\Pair_A(\CCC)$ has a $t$-structure such that
\[ \Pair_A(\CCC)_{\geq 0} \simeq \Pair_A(\CCC) \times_{\Pair(\CCC)} \Pair(\CCC)_{\geq 0}, \]
where $\Pair_A(\CCC) \to \Pair(\CCC)$ is the forgetful functor. This $t$-structure has the property that the adjunction $\LSym_{(A,0)} \dashv U_{(A,0)}$ restrict to an adjunction
\[ (\LSym_{(A,0)} \dashv U_{(A,0)}) \colon \Pair_A(\CCC)_{\geq 0} \rightleftarrows \DPair_A(\CCC)_{\geq 0} \]
where $\DPair_A(\CCC)_{\geq 0}$ by definition is the subcategory of $\DPair_A(\CCC)$ spanned by those $(A,0)$-algebras for which the underlying $(A,0)$-module is connective.

Likewise, canonically $(\id_A \colon A \to A) \in \DAlg(\Ar(\CCC))_{\geq 0}$, and we can relativize the derived algebraic context $\Ar(\CCC)$ to $(\id_A \colon A \to A)$. In this case we write
\[ \Ar_A(\CCC) \coloneqq \Mod_{(\id_A \colon A \to A)}(\Ar(\CCC)). \]

\begin{Prop}
    Taking cofibers induces a $t$-exact, symmetric monoidal equivalence
    \[ \Pair_A(\CCC) \to \Ar_A(\CCC) \colon (B,J) \mapsto (B \to B/J), \]
    which commutes with base-change in the sense that the diagram
    \begin{center}
        \begin{tikzcd}
            \Pair_A(\CCC) \arrow[r] \arrow[d] & \Ar_A(\CCC) \arrow[d] \\
            \Pair_{A'}(\CCC) \arrow[r] & \Ar_{A'}(\CCC)
        \end{tikzcd}
    \end{center}
    commutes, for any $A \to A'$ in $\DAlg(\CCC)_{\geq 0}$.
\end{Prop}

\begin{proof}
   The first statement is immediate by construction. The second statement follows from the fact that $(-) \otimes_A A'$ is exact. 
\end{proof}

\section{Quasi-coherent ideals \& Rees algebras}
We globalize the previous section to quasi-coherent ideals over a $\CCC_{\geq 0}$-stack $X$, for $\CCC$ the given derived algebraic context. We verify the familiar description of $B\G_m$ and $\Theta$ as classifying stacks for graded and filtered objects in \S\ref{Subsec:fil_grad_stacky}, introduce vanishing loci of quasi-coherent ideal in \S\ref{Subsec:V}, and apply our theory to the construction of scheme-theoretic images of arbitrary maps of $\CCC_{\geq 0}$-stacks in \S\ref{Subsec:ST_image}. A main result of this section is Corollary~\ref{Cor:affine_def}, where we show that the deformation space of $W \to X$ is almost affine whenever $W \to X$ is almost affine.
\subsection{Relative spectra and global sections}
Let $X$ be a stack. We have an adjunction
\[ \Gamma \dashv \widetilde{\Spec}_X(-) \colon \St_X \rightleftarrows \QAlg(X)^\op \]
where $\Gamma$ sends $f \colon W \to X$ to $f_*\OOO_W$. 
We obtain a factorization
\[ W \to \widetilde{\Spec}_X f_*\OOO_W \to X. \]

\begin{Rem}
\label{Rem:BC_rel_specnc}
    Let $g \colon X' \to X$ be a morphism of $\CCC_{\geq 0}$-stacks and $\BB \in \QAlg(X)$. Then the canonical square
    \begin{center}
    \begin{tikzcd}
        \widetilde{\Spec}(g^*\BB) \arrow[r] \arrow[d] & \widetilde{\Spec}(\BB) \arrow[d] \\
        X' \arrow[r, "g"] & X
    \end{tikzcd}
    \end{center}
    is Cartesian. This is immediate from comparing universal properties.
\end{Rem}

\begin{Lem}
\label{Lem:affine}
    Consider a Cartesian diagram
    \begin{center}
        \begin{tikzcd}
            W' \arrow[r, "g"] \arrow[d, "p"] & X' \arrow[d, "q"] \\
            W \arrow[r, "f"] & X
        \end{tikzcd}
    \end{center}
    where 
    \begin{enumerate}
        \item $p,q$ are effective epimorphisms,
        \item the base-change transformation
    $q^*f_* \to g_*p^* $
    is an equivalence,
    \item $f_*$ is $t$-exact, and
    \item $g$ is affine.
    \end{enumerate}
Then, $f$ is affine.
\end{Lem}

\begin{proof}
		By base-change it holds
		\[ q^*f_* \OOO_W \simeq g_*p^*\OOO_W \simeq g_* \OOO_{W'} \simeq \tau_{\geq 0}g_* \OOO_{W'}. \]
		Since $f_*$ is $t$-exact, we thus obtain
		the Cartesian diagram
		\begin{center}
			\begin{tikzcd}
				W' \arrow[r, "g'"] \arrow[d, "p"] & \Spec_{X'} (g_*\OOO_{W'}) \arrow[d, "r'"] \arrow[r]  & X' \arrow[d, "q"] \\
				W \arrow[r, "f'"] & \Spec_X (f_*\OOO_W) \arrow[r] & X,
			\end{tikzcd}
		\end{center}   
		where $g'$ is an equivalence and $p,r'$ are effective epimorphisms. It follows that $f'$ is an equivalence, hence that $f$ is affine.
\end{proof}

Write $\rho \colon * \to B\G_m$ for the universal $\G_m$-torsor, and $\pi \colon B\G_m \to *$ for the structure map. Then the adjunction $\rho^* \dashv \rho_*$ corresponds to the adjunction
\[ \mathrm{forget} \dashv (-) \otimes \mathrm{1}_\CCC[t,t^{-1}] \colon \CCC^\Z \rightleftarrows \CCC, \]
and $\pi^* \dashv \pi_*$ corresponds to the adjunction
\[ \mathrm{trivial} \dashv (-)_0 \colon \CCC \rightleftarrows \CCC^\Z. \]

\begin{Cor}
    The morphism $\rho \colon * \to B\G_m$ is affine, as is $i \colon B\G_m \to \Theta$.
\end{Cor}

\begin{proof}
    We have a Cartesian diagram
    \begin{center}
        \begin{tikzcd}
            \G_m \arrow[r, "p"] \arrow[d, "p"] & * \arrow[d, "\rho"] \\
            * \arrow[r, "\rho"] & B\G_m
        \end{tikzcd}
    \end{center}
    From the description of $\QCoh(B\G_m)$ as $\CCC^\Z$, it follows that $\rho^*\rho_* \simeq p_*p^*$, and that $\rho_*$ is $t$-exact. Hence $\rho$ is affine by Lemma \ref{Lem:affine}.

    For the second statement, consider the Cartesian diagram
    \begin{center}
        \begin{tikzcd}
            \{0\} \arrow[r] \arrow[d, "\rho"] & \A^1 \arrow[d, "q"] \\
            B\G_m \arrow[r, "i"] & \Theta.
        \end{tikzcd}
    \end{center}
    Under the equivalences of stable symmetric monoidal categories with $t$-structures $\QCoh(B\G_m) \simeq \CCC^\Z$ and $\QCoh(\Theta) \simeq \Mod_{\mathrm{1}_\CCC[t^{-1}]}(\CCC^\Z)$, it holds that $i_*$ is the $\Z$-graded restriction of scalars functor along $\mathrm{1}_\CCC[t^{-1}] \to \mathrm{1}_\CCC$, and $q^*$ is the functor which forgets the grading. We are thus again in a situation where Lemma \ref{Lem:affine} applies, and the claim follows.
\end{proof}
By definition of virtual Cartier divisors, we conclude:
\begin{Cor}
    Any virtual Cartier divisor $D \to T$ is affine.
\end{Cor}

\subsection{Filtered and graded objects via classifying stacks}
\label{Subsec:fil_grad_stacky}
 For $B \in \DAlg(\CCC^\Z)$, we write $\Mod^\Z_B$ for the category of $\Z$-graded $B$-modules in $\CCC$.
 
For the filtered counterpart of the equivalence $\CCC^\Z \simeq \QCoh(B\G_{m})$, consider the \emph{Rees functor}
 \[ \varphi \colon \Mod^\Z_{\mathrm{1}_\CCC[t^{-1}]} \to \Fil(\CCC) \]
 which sends $N \in \Mod^\Z_{\mathrm{1}_\CCC[t^{-1}]}$ to the fitered object $\varphi N$ such that $\ev^i(\varphi N) \coloneqq N^i$, and with structure morphisms
 \[ \times t^{i' - i } \colon \ev^{i}(\varphi N) \to \ev^{i'}(\varphi N) \]
 for $i' \leq i$.

The following goes back to Simpson.
 \begin{Prop}
 \label{Prop:Fil_is_Ztalg}
     The functor $\varphi \colon \Mod^\Z_{\mathrm{1}_\CCC[t^{-1}]} \to \Fil(\CCC) $ is an equivalence of derived algebraic contexts, and thus exhibits $\Theta$ as the classifying stack for filtered objects in $\CCC$ via the equivalences
     \[ \QCoh(\Theta) \simeq \Mod^\Z_{\mathrm{1}_\CCC[t^{-1}]} \simeq_\varphi \Fil(\CCC). \]
 \end{Prop}

\begin{proof}
    We have already seen that $\QCoh(\Theta) \simeq \Mod^\Z_{\mathrm{1}_\CCC[t^{-1}]}$, and the second equivalence is \cite[Prop.~3.2.9]{RaksitDerived}.
\end{proof}

\begin{Cor}
\label{Cor:gr_via_pullback}
    Under the equivalences $\Fil(\CCC) \simeq \QCoh(\Theta)$ and $\CCC^\Z \simeq \QCoh(B\G_m)$, the functor 
    \[ \gr \colon \Fil(\CCC) \to \CCC^\Z \]
    corresponds to the functor
    \[ i^* \colon \QCoh(B \G_m) \to \QCoh(\Theta) \]
    induced by pulling back along the closed immersion $i \colon B \G_m \to \Theta$. Likewise, the forgetful functor (i.e., the colimit functor)
     \[ U \colon \Fil(\CCC) \to \CCC \]
    corresponds to the functor
    \[ j^* \colon \QCoh(\Theta) \to \CCC \]
    induced by the pullback along $j \colon * \to \Theta$.   
\end{Cor}

\begin{proof}
    The first statement is \cite[Prop.\ 3.2.9]{RaksitDerived}. For the second statement it suffices to observe that the pushforward map
    \[ j_* \colon \Mod^\Z_{\mathrm{1}_\CCC[t,t^{-1}]} \to \Mod^\Z_{{\mathrm{1}_\CCC}[t^{-1}]} \]
    corresponds to the diagonal functor $\CCC \to \Fil(\CCC)$ under the equivalences $\Mod^\Z_{\mathrm{1}_\CCC[t,t^{-1}]} \simeq \QCoh([\G_m/\G_m]) \simeq \CCC$ and $\Mod^\Z_{\mathrm{1}_\CCC[t^{-1}]} \simeq \QCoh(\Theta) \simeq \Fil(\CCC)$.
\end{proof}

\begin{Rem}
    Note that $\QAlg(\Theta) \simeq \DAlg^\Z_{{\mathrm{1}_\CCC}[t^{-1}]}$, hence also $\QAlg(\Theta) \simeq \DAlg(\Fil(\CCC))$. 
\end{Rem}

\subsection{Deformation to the normal bundle via Rees algebras}
\label{Subsec:DeftoNB}

\begin{Lem}
\label{Lem:D2NB_Rees}
    Let $A \to B$ be a morphism of derived $\CCC$-algebras, with $A$ connective. Then there is a canonical $\G_{m}$-equivariant equivalence
    \[ D_{\widetilde{\Spec} B/\Spec A} \simeq \widetilde{\Spec} R^{\ext}_{B/A} \]
    In particular, if $A,B$ are  both connective, then the extended Rees algebra of $A \to B$ as defined in \cite{OrenJeroenblowups} coincides with our definition.
\end{Lem}

\begin{proof}
    It suffices to show that the composition
    \begin{multline*}
        G \colon \DAlg^\Z_{\mathrm{1}_\CCC[t^{-1}]}\simeq \DAlg(\Fil(\CCC)) \xrightarrow{(-)^{\geq 0}} \\
        \DAlg(\Fil^{\geq 0}(\CCC)) \xrightarrow{\res} \DPair(\CCC) \simeq \DAlg(\Ar(\CCC))
    \end{multline*}
    is equivalent to the functor that sends $Q \in \DAlg^\Z_{\mathrm{1}_\CCC[t^{-1}]}$ to
    \[ (Q^0 \to  (Q/(t^{-1}))^0 )\coloneqq (\ev^0 Q \to \ev^0(Q \otimes_{\mathrm{1}_\CCC[t^{-1}]} \mathrm{1}_\CCC)) \]
    And indeed this is so, as follows by tracing through the definitions.
\end{proof}

\begin{Cor}
    For any $A \in \DAlg(\CCC)_{\geq 0}$, the functor $B \mapsto R^\ext_{B/A}$ is a fully faithful left adjoint to the functor $\DAlg^\Z_{A[t^{-1}]} \to \DAlg_A \colon Q \mapsto (Q/(t^{-1}))^0$.
\end{Cor}

\begin{proof}
    It is clear that $R^{\ext}_{(-)/A}$ is left adjoint to $Q \mapsto (Q/(t^{-1}))^0$. Fully faithfulness follows from Lemma \ref{Lem:Rees_fully_faithful}, together with Proposition \ref{Prop:Fil_is_Ztalg}. 
\end{proof}

\begin{Cor}
\label{Cor:NB_is_NC}
    Let $A \to B$ be a surjective morphism of $\CCC_{\geq 0}$-algebras, with ideal $I$. Then there is a canonical equivalence
    \[ \gr(R^\ext_{B/A}) \simeq \LSym_{B}(L_{B/A}[-1](-1)). \]
    In particular, it holds that $I/I^2 \simeq \gr^1(R^\ext_{B/A}) \simeq L_{B/A}[-1]$.
\end{Cor}

\begin{proof}
    This follows from Lemma~\ref{Lem:D2NB_Rees} by Corollary~\ref{Cor:gr_via_pullback}, together with base-change for Cartesian squares of affine morphisms and fully faithfulness of $\Spec(-)$ in the connective case.
\end{proof}

\begin{Def}
    Let $f \colon W \to X$ be a morphism of $\CCC_{\geq 0}$-stacks and write $\pi \colon \cD_{W/X} \to X \times \Theta$ for the structure map. Then the \emph{extended Rees algebra} of $f$ is the $\Z$-graded, quasi-coherent $\OOO_X[t^{-1}]$-algebra
    \[ \RR_{W/X}^\ext \coloneqq \pi_* \OOO_{\cD_{W/X}} \in \QAlg(X \times \Theta) \simeq \QAlg^\Z(X)_{\OOO_X[t^{-1}]/}. \]
\end{Def}

\begin{Cor}
\label{Cor:affine_def}
    Let $\BB \in \QAlg(X)$ and put $W \coloneqq \widetilde{\Spec}_X(\BB)$. 
    If $X$ is affine, then canonically $\RR^\ext_{W/X} \simeq R_{\OO_W/\OO_X}^\ext$ under $\QCoh(X) \simeq \Mod_{\OO_X}(\CCC)$.
    In general, the canonical $\G_m$-equivariant map
    \[ D_{W/X} \to \widetilde{\Spec}_{X \times \A^1} \RR^\ext_{W/X} \]
    is an equivalence. Moreover, if $W \to X$ is a closed immersion, then $D_{W/X} \to X \times \A^1$ is affine. 
\end{Cor}

\begin{proof}
    The first and second claim follow from Lemma~\ref{Lem:D2NB_Rees} and Remark~\ref{Rem:BC_rel_specnc}. The third claim now is \cite[Thm.~4.22]{OrenJeroenblowups}.
\end{proof}

\subsection{Quasi-coherent pairs}
Let $\PrL_\st$ be the full subcategory of $\PrL$ spanned by presentable categories which are stable. Observe that $\PrL_\st$ is a symmetric monoidal subcategory of $\PrL$ (with respect to the Lurie tensor product).
\begin{Def}
    Write $\Pair \colon \St^\op \to \CAlg(\PrL_\st)$ for the right Kan extension along the Yoneda embedding of the functor $\Aff^\op \to \CAlg(\PrL_\st)$ which sends $\Spec A$ to $\Pair_A(\CCC)$.
\end{Def}
For a stack $X$ write $\Pair(X)_{\geq 0}$ for the full subcategory of $\Pair(X)$ spanned by those $\FF \in \Pair(X)$ such that $\varphi^* \FF \in \Pair_A(\CCC)_{\geq 0}$ for all $\varphi \colon \Spec(A) \to X$. 
\begin{Prop}
    The category $\Pair(X)_{\geq 0} \subset \Pair(X)$ is the connective part of a $t$-structure on $\Pair(X)$.
\end{Prop}
\begin{proof}
    The canonical diagram
    \begin{center}
        \begin{tikzcd}
            \Pair(X)_{\geq 0} \arrow[r] \arrow[d] & \Pair(X) \arrow[d] \\
            \QCoh(X)_{\geq 0} \times \QCoh(X)_{\geq 0} \arrow[r] & \QCoh(X) \times \QCoh(X)
        \end{tikzcd}
    \end{center}
    is Cartesian, where the vertical arrows are the forgetful functors $\FF \mapsto (\FF^0,\FF^1)$. Indeed the condition holds affine-locally by definition, and is stable under base-change hence survives the globalization via right Kan extension. It follows that $\Pair(X)_{\geq 0} \subset \Pair(X)$ is closed under extensions, hence forms the connective part of a $t$-structure on $\Pair(X)$ \cite[Prop.~1.2.1.16, Prop.~1.4.4.11]{LurieHA}.
\end{proof}
\begin{Rem}
    Using the pointwise formula for right Kan extensions, one can show that
    \[ \Pair(X) \simeq \Fun(([0,1],\leq,+)^\op,\QCoh(X)), \]
    where the right-hand side is endowed with the neutral $t$-structure and the Day convolution symmetric monoidal structure from \S \ref{subsec:first_defs}.
\end{Rem}
\begin{Def}
    Write $\DPair \colon \St^\op \to \CAlg(\PrL)$ for the right Kan extension along the Yoneda embedding of the functor $\Aff^\op \to \CAlg(\PrL)$ which sends $\Spec A$ to $\DPair_A(\CCC)$.
\end{Def}
The monadic adjuncion $(\LSym_{(A,0)} \dashv U_{(A,0)}) \colon \Pair_A(\CCC) \rightleftarrows \DPair_A(\CCC)$ commutes with base-change along maps of connective $\CCC$-algebras $A \to A'$ in the obvious way. We thus obtain a canonical monadic adjunction
\[ (\LSym_{(\OOO_X,0)} \dashv U_{(\OOO_X,0)}) \colon \Pair(X) \rightleftarrows \DPair(X) \]
for any $\CCC_{\geq 0}$-stack $X$.

Write $\DPair(X)_{\geq 0} \coloneqq \DPair(X) \times_{\Pair(X)} \Pair(X)_{\geq 0}$. Then the inclusion $\DPair(X)_{\geq 0} \to \DPair(X)$ has a right adjoint
\[ \tau_{\geq 0} \colon \DPair(X) \to \DPair(X)_{\geq 0} \]
which commutes with the forgetful functors and the connective cover functor $\tau_{\geq 0} \colon \Pair(X) \to \Pair(X)_{\geq 0}$ in the obvious way.

\begin{Lem}
\label{Lem:RKE_Ar}
The right Kan extension of the functor $\Spec(A) \mapsto \Ar_A(\CCC)$ along the Yoneda embedding is the functor
\[ \St^\op \to \CAlg(\PrL_\st) \colon X \mapsto \Arr(\QCoh(X)). \]
Likewise, the right Kan extension of $\Spec(A) \mapsto \DAlg_{(\id_A \colon A \to A)}(\Ar(\CCC))$ along the Yoneda embedding is the functor
\[ \St^\op \to \CAlg(\PrL) \colon X \mapsto \Arr(\QAlg(X)). \]
\end{Lem}

\begin{proof}
    This is immediate from the fact that $\Ar_A(\CCC) \simeq \Arr(\Mod_A(\CCC))$ and $\DAlg_{(\id_A \colon A \to A)}(\Ar(\CCC)) \simeq \Arr(\DAlg_A(\CCC))$.
\end{proof}

As before, the category $\Arr(\QCoh(X))$ has a $t$-structure such that the connective part  $\Arr(\QCoh(X))_{\geq 0}$  is determined section-wise. Write 
\[ \Sur(\QCoh(X)) \coloneqq \Arr(\QCoh(X))_{\geq 0}, \]
and put $\Sur(\QAlg(X)) \coloneqq \Arr(\QAlg(X)) \times_{\Arr(\QCoh(X))} \Sur(\QCoh(X))$.

\begin{Prop}
\label{Prop:Pair_Arr_global}
    Taking cofibers induces equivalences of categories
    \[ \Pair(X) \to \Arr(\QCoh(X)), \quad \textrm{and} \quad \DPair(X) \to \Arr(\QAlg(X)), \]
    which restrict to equivalences of categories
    \[ \Pair(X)_{\geq 0} \to \Sur(\QCoh(X)), \quad \textrm{and} \quad \DPair(X)_{\geq 0} \to \Sur(\QAlg(X)). \]    
\end{Prop}

\begin{proof}
    By Lemma~\ref{Lem:RKE_Ar}, the first statement follows from the discussion in \S \ref{subsec:Pair_base}. For the second statement, note that the inclusion $\Pair_A(\CCC)_{\geq 0} \to \Pair_A(\CCC)$ commutes with base-change along maps of connective $\CCC$-algebras $A \to A'$ in the obvious way, as does the inclusion $\Ar_A(\CCC)_{\geq 0} \to \Ar_A(\CCC)$. Hence the corresponding affine statement survives to the global statement via right Kan extension.
\end{proof}

A morphism of connective objects is called \emph{surjective} if it is surjective on $\pi_0$. From Proposition~\ref{Prop:Pair_Arr_global} it is immediate that $\Sur(\QCoh(X))$ is spanned by the surjective morphisms of connective quasi-coherent $\OOO_X$-modules. Likewise for $\Sur(\QAlg(X))$.

\subsection{Quasi-coherent ideals}
Let $X$ be a $\CCC_{\geq 0}$-stack. An object of $\DPair(X)$ can be described as a pair $(\BB,\JJ)$ where $\BB$ is a quasi-coherent $\OOO_X$-algebras and $\JJ \to \BB$ is a quasi-coherent ideal of $\BB$ in the sense that $(f^*\BB,f^*\JJ)$ is an ideal pair in $\DPair_A(\CCC)$ for any $f \colon \Spec A \to X$. 
\begin{Def}
    For a stack $X$ let $\Ideal(X)$ be the full subcategory of $\DPair(X)$ spanned by objects of the form $(\OOO_X,\II)$. Objects of $\Ideal(X)$ are called \emph{quasi-coherent ideals (on $X$)}, and are simply written $\II$. 
\end{Def}
The category $\Ideal(X)$ is a symmetric monoidal subcategory of $\DPair(X)$. This is immediate from the description of the symmetric monoidal structure on $\DPair(X)$ as the global version of the pushout product symmetric monoidal structure. Moreover, writing $\Ideal(X)_{\geq 0} \coloneqq \Ideal(X) \times_{\DPair(X)} \DPair(X)_{\geq 0}$, the inclusion $\Ideal(X)_{\geq 0} \to \Ideal(X)$ has a right adjoint
\[ \tau_{\geq 0} \colon \Ideal(X) \to \Ideal(X)_{\geq 0}.\]
Note also that there is a canonical forgetful functor $\Ideal(X) \to \QCoh(X)$ which restricts to $\Ideal(X)_{\geq 0} \to \QCoh(X)_{\geq 0}$ and commutes with $\tau_{\geq 0}$ in the obvious way. 
\begin{Lem}
\label{Lem:Idea_QAlg}
    The equivalence $\DPair(X) \xrightarrow{\sim} \Arr(\QAlg(X))$ induces a symmetric monoidal equivalence
    \[ \Ideal(X) \to \QAlg(X) \colon \II \mapsto \OOO_X/\II \]
    such that the underlying quasi-coherent $\OOO_X$-module of $\OOO_X/\II$ is the cofiber of the underlying map $\II \to \OOO_X$ in $\QCoh(X)$. 
\end{Lem}
\begin{proof}
    This is immediate from Proposition~\ref{Prop:Pair_Arr_global}.
\end{proof}

\begin{Warn}
    As is to be expected from Lemma~\ref{Lem:Idea_QAlg}, the category $\Ideal(X)$ is not stable in general. Even sillily so: if $0$ is terminal in $\Ideal(X)$ then $X = \emptyset$.
\end{Warn}

\subsection{Vanishing loci}
\label{Subsec:V}
Let $X$ be a $\CCC_{\geq 0}$-stack. The adjunction 
\[ \Gamma \dashv \widetilde{\Spec}_X(-) \colon \St_X \rightleftarrows \QAlg(X)^\op\]
induces an equivalence between $\QAlg(X)^\op_{\geq 0}$ and the full subcategory $\Aff(X) \subset \St_X$ of stacks affine over $X$. 
\begin{Def}
\label{Def:almost_affine}
   Write $\widetilde{\Aff}(X)$ for the full subcategory of $\St_X$ spanned by $X$-stacks 
of the form $\widetilde{\Spec}_X(\AA) \to X$. Let us call $T \to X$ \emph{almost affine} if it lives in $\widetilde{\Aff}(X)$. 
\end{Def}

\begin{Warn}
\label{Warn:nc_affine}
    The functor $\widetilde{\Spec}_X(-)$ in general is \emph{not} fully faithful. Morally, this is because we are probing nonconnective geometry using `only' connective $\CCC$-algebras (hence this can also go wrong for $X=*$). 
\end{Warn}    

\begin{Rem}   
\label{Rem:affine_stacks}
    Notwithstanding Warning~\ref{Warn:nc_affine} and as explained in \cite[Ex.~3.8]{OrenJeroenblowups}, there is always a minimal choice of extending the given topology on $\DAlg(\CCC)_{\geq 0}^\op$ to a subcanonical topology on $\DAlg(\CCC)^\op$ for which $\Mod_{(-)}$ still satisfies descent, and such that we obtain a fully faithful functor
    \[ \St(\CCC) \triangleq \Sh(\DAlg(\CCC)_{\geq 0}^\op) \xhookrightarrow{i} \Sh(\DAlg(\CCC)^\op) \]
    for the categories of sheaves with respect to these topologies, together with adjunctions $t \dashv i \dashv \rho$.\footnote{The notation of \cite{OrenJeroenblowups} is $\St_{\CCC_{\geq 0}} \coloneqq \Sh(\DAlg(\CCC)_{\geq 0}^\op)$ and $\St_{\CCC} \coloneqq \Sh(\DAlg(\CCC)^\op$, similarly for affine objects.} The functor $t$ preserves affines (that is, \emph{nonconnective} affines now), on which it is given by taking connective covers on global sections. The functor $\rho$ is induced by restricted Yoneda. 
    By design, we have a commutative (even Cartesian) diagram
    \begin{center}
        \begin{tikzcd}
            \DAlg(\CCC)_{\geq 0}^\op \arrow[r, hookrightarrow] \arrow[d, hookrightarrow, "\Spec"] & \DAlg(\CCC)^\op \arrow[d, hookrightarrow, "\Spec^\nc"] \\
            \St(\CCC) \arrow[r, hookrightarrow, "i"] & \Sh(\DAlg(\CCC)^\op),
        \end{tikzcd}
    \end{center}
    where $\Spec^\nc$ is induced by Yoneda, and a factorization
    \begin{center}
        \begin{tikzcd}
             & \DAlg(\CCC)^\op \arrow[d, hookrightarrow, "\Spec^\nc"] \arrow[dl, "\widetilde{\Spec}", swap] \\
            \St(\CCC)  & \Sh(\DAlg(\CCC)^\op) \arrow[l, "\rho"].
        \end{tikzcd}
    \end{center}
    
    Since $i$ is fully faithful and $i \dashv \rho$, the failure of $\widetilde{\Spec}$ being fully faithful seems as deep as it gets. Remarkably though, when restricted to derived $\CCC$-algebras which are bounded (either from above or below), then $\widetilde{\Spec}$ is fully faithful (over $\Mod_\Z$ and with the fpqc topology)---see \cite[Thm.~3.4]{MathewAffine}, which goes back to \cite[Prop.~2.2.2]{ToenChamps}. Moreover, when restricted to $n$-truncated derived rings, the essential image of $\widetilde{\Spec}$ consists  exactly of the $n$-derived, affine stacks---see \cite[Cor.~3.6]{MathewAffine} and references therein. Here, `affine stacks' refers to the terminology introduced in \cite{ToenChamps}, and a main point from \cite{MathewAffine} is the application of the theory of derived rings to the idea of affine stacks.

    Following \cite[Def.~3.24]{OrenJeroenblowups}, we call a morphism $f' \colon W' \to X'$ in the category $\Sh(\DAlg(\CCC)^\op)$ \emph{nonconnectively affine} if $W' \simeq \Spec^\nc_{X'}(f_*\OO_{W'})$, and we write $\Aff^\nc(X)$ for the category of such. If $X' = iX$ for $X \in \St(\CCC)$, then we have induced essentially surjective functors
    \[ \Aff^\nc(X') \xrightarrow{\rho} \widetilde{\Aff}(X) \xrightarrow{t} \Aff(X), \]
    which in general are not fully faithful. In other words, for given $f \colon W \to X$ in $\St(\CCC)$ it holds: if $if$ is nonconnectively affine then $f$ is almost affine, but the converse fails in general.
\end{Rem} 

\begin{Rem}
\label{Rem:nc_affine_D}
    We have seen in Corollary~\ref{Cor:affine_def} that the deformation space $D_{W/X}$ is almost affine over the base $X$ whenever $W \to X$ is an almost affine morphism of $\CCC_{\geq 0}$-stacks. With the same argument but now carried out in the category $\Sh(\DAlg(\CCC)^\op)$ from Remark~\ref{Rem:affine_stacks}, one shows that $D_{W'/X'}$ is nonconnectively affine over $X'$ for any nonconnectively affine map $W' \to X'$ in $\Sh(\DAlg(\CCC)^\op)$.
\end{Rem}

\begin{Def}
    The \emph{vanishing locus} of a quasi-coherent ideal $\II$ on $X$ is the stack
    \[ V(\II) \coloneqq {\widetilde{\Spec}_X(\OOO_X/\II)}. \]
    This exhibits a functor $V(-) \colon \Ideal(X) \to \widetilde{\Aff}(X)$.
\end{Def}
\begin{Rem}
\label{Rem:sums_of_ideals}
    The symmetric monoidal structure on $\Pair(X)$ induces the coCartesian symmetric monoidal structure on $\DPair(X)$. The restriction to $\Ideal(X)$ is written as $+$. For $\II,\JJ \in \Ideal(X)$ it holds that
    \[ V(\II + \JJ) \simeq V(\II) \times_X V(\JJ) \]
    by Lemma~\ref{Lem:Idea_QAlg}, making $V(-)$ symmetric monoidal.
\end{Rem}
Recall that a \emph{closed immersion} is an affine morphism $Z \to X$ such that $\OOO_X \to \OOO_Z$ is surjective.
\begin{Lem}
    Let $\II$ be a quasi-coherent ideal on $X$. Then the canonical map $f \colon V(\II) \to X$ is
    \begin{enumerate}
        \item affine if and only if $\II$ is $(-1)$-connective, and
        \item a closed immersion if and only if $\II$ is connective.
    \end{enumerate}
\end{Lem}
\begin{proof}
    This follows immediately from the exact sequence $\II \to \OOO_X \to \OOO_X/\II$ in $\QCoh(X)$.
\end{proof}
For any map $W \to X$ we can associate to it a quasi-coherent ideal $\II(W)$ via the composition
\[ \St(X) \xrightarrow{\Gamma} \QAlg(X)^\op \simeq \Ideal(X)^\op. \]

Write $\Cld(X)$ for the full subcategory of $\Aff(X)$ spanned by closed immersions $Z \to X$. 
\begin{Cor}
\label{Cor:V_Cld_equiv}
    Taking vanishing loci exhibits an equivalence
    \[ \Ideal(X)_{\geq 0} \to \Cld(X) \colon \II \mapsto V(\II) \]
    between connective quasi-coherent ideals and closed substacks of $X$.\footnote{Here and in what follows, a \emph{closed substack of $X$} is by definition a closed immersion $Z \to X$, even though in general a closed substack is not a subobject.} 
\end{Cor}

\begin{Rem}
\label{Rem:inf_nbhd}
    Let $(A,I)$ be an ideal pair and $A \to A'$ a morphism in $\DAlg(\CCC)$. Write $(A',I')$ for the ideal pair obtained by base-change of $A \to A/I$ along $A \to A'$. Since the Rees algebra is stable under arbitrary pullback, we obtain a canonical equivalence
    \[ (A/I^{n}) \otimes_A A' \simeq A'/(I')^{n} \]
    In particular, the construction $(A,I) \mapsto (A,I^{n})$ globalizes in the obvious way to a functor
    \[ \Ideal(X) \to \Ideal(X) \colon \II \mapsto \II^n, \]
    for any $X \in \St$. If $Z \coloneqq V(\II)$ is the vanishing locus of $\II$, then we write
    \[ Z^{(n)} \coloneqq V(\II^{n+1}) \to X \]
    for the vanishing ideal of the corresponding ideal power, and call it the \emph{$n$th infinitesimal neighborhood} of $Z$ in $X$. 
    
    If $Z \to X$ is a closed immersion, then we can concretely describe $\OOO_{Z^{(n)}}$ as
    \[ \OOO_{Z^{(n)}} \simeq (\RR^\ext(\OOO_X,\II)/(t^{-n-1}))_0. \]
    In particular, for any map $X' \to X$ we obtain a Cartesian diagram
    \begin{center}
        \begin{tikzcd}
            (Z \times_X X')^{(n)} \arrow[d] \arrow[r] & X' \arrow[d] \\
            Z^{(n)} \arrow[r] & X,
        \end{tikzcd}
    \end{center}
    again by stability of the Rees algebra under abitrary pullback.
\end{Rem}

\begin{Cor}
\label{Cor:cotangent_Quillen}
    Let $Z \to X$ be a closed immersion of $\CCC_{\geq 0}$-stacks, with corresponding ideal $\II \in \Ideal(X)$. Then there is a canonical equivalence
    \[ \LL_{Z/X}[-1] \simeq \II/\II^2. \]
\end{Cor}

\begin{proof}
    Since the deformation to the normal bundle is stable under base change, this reduces to the affine case, which is Corollary~\ref{Cor:NB_is_NC}.
\end{proof}

\subsection{Scheme-theoretic images}
\label{Subsec:ST_image}
Let $X$ be a $\CCC_{\geq 0}$-stack. By Proposition~\ref{Prop:Pair_Arr_global} the inclusion $\Sur(\QAlg(X)) \to \Arr(\QAlg(X))$ has a right adjoint. This right adjoint sends a given morphism $\varphi \colon \AA \to \BB$ of quasi-coherent $\OOO_X$-algebras to the morphism $\tau_{\geq 0} \AA \to \varphi \BB$ of quasi-coherent $\OOO_X$-algebras which fits into a commutative diagram  with exact rows
	\begin{equation}
		\label{Eq:Image}
		\begin{tikzcd}
			\tau_{\geq 0} \II \arrow[d] \arrow[r] & \tau_{\geq 0} \AA \arrow[r] \arrow[d] & \varphi \AA \arrow[d] \\
			\II \arrow[r] & \AA \arrow[r, "\varphi"] & \BB
		\end{tikzcd}
	\end{equation}
 in $\QCoh(X)$, where the square on the left is induced by taking connective covers in $\QCoh(X)$, and the square on the right lives in $\QAlg(X)$. We call $\varphi \AA$ the \emph{image} of $\varphi$. 
\begin{Lem}
\label{Lem:image_ring_map}
    Let $\varphi \colon \AA \to \BB$ be a morphism of quasi-coherent $\OOO_X$-algebras with ideal $\II$, and suppose that $\AA$ is connective. Let $\JJ$ be the fiber of the map $\varphi \AA \to \BB$. Then:
    \begin{enumerate}
        \item there is a canonical equivalence $\JJ \simeq \tau_{\leq -1} \II$,
        \item $\varphi \AA$ is connective and $\AA \to \varphi \AA$ is surjective, and
        \item if $\BB$ is also connective then $\varphi \AA\to \BB$ is injective on $\pi_0$ and an equivalence on $\pi_n$ for $n > 0$.
    \end{enumerate}
\end{Lem}
\begin{proof}
    All three statements follow from the exact diagram
    \begin{center}
        \begin{tikzcd}
            \tau_{\geq 0} \II \arrow[r] \arrow[d] & \II \arrow[d] \arrow[r] & \AA \arrow[d] \\
            0 \arrow[r] & \JJ \arrow[d] \arrow[r] & \varphi \AA \arrow[d] \\
            & 0 \arrow[r] & \BB. 
        \end{tikzcd}
    \end{center} 
\end{proof}
\begin{Def}
    For a morphism $f \colon W \to X$ of stacks define the \emph{scheme-theoretic image} of $f$ as the factorization
    \[ W \to \Spec_X(\varphi \OOO_X) \to X \]
    where $\varphi \colon \OOO_X \to f_* \OOO_W$ is the canonical map of quasi-coherent $\OOO_X$-algebras.
\end{Def}
We also write $f(W) \coloneqq \Spec_X(\varphi \OOO_X)$ for the scheme-theoretic image of $f \colon W \to X$ (suppressing the factorization from notation). Note that $f(W)$ is also the vanishing locus $V(\tau_{\geq 0}\II)$ of the connective cover of the quasi-coherent ideal $\II$ associated to the map  $\widetilde{\Spec}_X(f_* \OOO_W) \to X$ induced by $f$. In particular, $f(W) \to X$ is a closed immersion. In fact, it holds:

\begin{Prop}
\label{Prop:scheme_theoretic_image}
    Let $f \colon W \to X$ be a morphism of stacks with scheme-theoretic image $i \colon Z \to X$. Then $Z$ is the smallest closed substack of $X$ through which $f$ factors in the sense that it exhibits an initial object in the category of factorizations of $f$ of the form
    \[ W \xrightarrow{f'} Z' \xrightarrow{i'} X \]
    where $i'$ is a closed immersion.
\end{Prop}

\begin{proof}
    Let $i' \colon Z' \to X$ be a closed immersion with quasi-coherent ideal $\II'$, and let $\II$ be the quasi-coherent ideal associated to $\widetilde{\Spec}_X(f_* \OOO_W) \to X$. Then we have natural equivalences
    \begin{align*}
        \St_X(W,Z') & \simeq \QAlg(X)(\OOO_X/\II', f_* \OOO_W) \\
        & \simeq \Ideal(X)(\II', \II) \\
        & \simeq \Ideal(X)_{\geq 0}(\II',\tau_{\geq 0}\II) \\
        & \simeq \St_X(V(\tau_{\geq 0}\II),V(\II')) \\
        & \simeq \St_X(Z,Z'),
    \end{align*}
    by the universal property of $\Spec_X$, the equivalence $\QAlg(X) \simeq \Ideal(X)$, the definition of the connective cover of quasi-coherent ideals, Corollary~\ref{Cor:V_Cld_equiv}, and the definition of the scheme-theoretic image. The claim follows.
\end{proof}

\begin{Warn}
    The scheme-theoretic image is \emph{not} obtained through the factorization system on $\St(\CCC)$ of effective epimorphisms followed by monomorphisms. Indeed, closed immersions are in general not monomorphisms in the derived settings. But even in the classical picture this goes wrong, since a non-empty open subscheme of an integral scheme is dense (in the algebraic setting). 
\end{Warn}

	\begin{Exm}
    \label{Ex:GR_image}
		Consider the context $\CCC = \Mod_k$, where $k$ is a discrete ring such that $\Q \subset k$. Then the scheme-theoretic image just defined coincides with the one found in \cite[\S 5.1]{GaitsgoryStudy}. Note that the assumption $\Q \subset k$ is necessary for the construction in \cite{GaitsgoryStudy} to make sense.
	\end{Exm}	

\begin{Exm}
    Write $\Aff^\heartsuit$ for the full subcategory of $\Aff$ spanned by affine schemes of the form $\Spec A$ where $A$ is discrete. Consider the adjunction
    \[ (j \dashv (-)_\cl) \colon \Aff^\heartsuit \rightleftarrows \Aff \]
    induced by $\pi_0$. Suppose that $\Aff^\heartsuit$ is endowed with a Grothendieck topology $J^\heartsuit$ such that both $j$ and $(-)_\cl$ preserve covering families (with respect to the fixed topology $J$ on $\Aff$). Write $\St^\cl$ for the full subcategory of $\PPP(\Aff^\heartsuit)$ spanned by the presheaves that satisfy $J^{\heartsuit}$-descent. Then we have an adjunction
    \[ (\iota \dashv (-)_\cl) \colon \St^\cl \rightleftarrows \St \]
    such that $(-)_\cl$ is given by precomposition with $j$, and $\iota$ is induced by left Kan extension and localization and is fully faithful. We call a $\CCC_{\geq 0}$-stack \emph{classical} if it is in the essential image of $\iota$.

    Let $f \colon W \to X$ be a morphism of stacks. Then the underlying classical stack $f(W)_\cl$ of the scheme-theoretic image of $f$ is the initial object in the category of factorizations
    \[ W_\cl \xrightarrow{f'} Z' \xrightarrow{i'} X_\cl \]
    of $f$ such that $i' \colon Z' \to X_\cl$ is a closed immersion and $Z'$ is classical. In particular, in the algebraic setting the underlying classical stack of the scheme-theoretic image recovers the classical definition of scheme-theoretic image.
\end{Exm}

\begin{Exm}
    Consider the algebraic setting $\CCC = \Mod_\Z$. Consider a morphism $A \to B$ of connective $\Z$-algebras. Write $\varphi \colon R^\ext_{B/A} \to A[t,t^{-1}]$ for the morphism obtained by inverting $t^{-1}$ (induced by the deformation to the normal bundle). Then the image $\varphi R_{B/A}^\ext \to A[t,t^{-1}]$ is spanned (as a space) by the classical extended Rees algebra of $\pi_0A \to \pi_0B$ in the sense that the canonical diagram
    \begin{center}
        \begin{tikzcd}
            \varphi R^{\ext}_{B/A} \arrow[r] \arrow[d] & A[t,t^{-1}] \arrow[d] \\
            \pi_0 A[I_0t,t^{-1}] \arrow[r] & \pi_0 A[t,t^{-1}]
        \end{tikzcd}
    \end{center}
    is Cartesian, where $I_0$ is the ideal of $\pi_0A \to \pi_0B$ in the classical sense. In particular, $\varphi R_{B/A}^\ext$ only sees the underlying classical map $\pi_0A \to \pi_0 B$. 
\end{Exm}

\section{$I$-complete notions}
We still let $\CCC$ be a derived algebraic context.
Throughout this section, we also let $X$ be a $\CCC_{\geq 0}$-stack and $\II$ a quasi-coherent ideal on $X$. We introduce $\II$-complete modules in \S\ref{Subsect:I_adic_complete}, $\II$-adically complete modules in \S\ref{Subsec:Derived_I_cmplt}, and show that these agree whenever $\II$ is locally finitely generated in Proposition~\ref{Prop:derived_I_complete_T}. We close with a definition of formal spectra and formal completions in \S\ref{Subsec:formal_spec}, which we again show to be the same whenever $\II$ is locally finitely generated in Theorem~\ref{Thm:Spf_is_completion}.
\subsection{$\II$-complete quasi-coherent modules}
\label{Subsect:I_adic_complete}
\begin{Def}
    A morphism $\NN' \to \NN$ of quasi-coherent $\OOO_X$-modules is called an \emph{$\OOO_X/\II$-equivalence} if $\NN' \otimes \II \to \NN \otimes \II$ is an equivalence. 
\end{Def}
Recall that we have an exact sequence $\II \to \OOO_X \to \OOO_X/\II$ in $\QCoh(X)$, and that $\OOO_X \to \OOO_X/\II$ naturally lives in $\QAlg(X)$. For $\MM \in \QCoh(X)$ we write $\II \MM \coloneqq \II \otimes \MM$, and let $\MM/\II\MM$ be the cofiber
\[ \II\MM \to \MM \to \MM/\II\MM. \]
Then $\NN' \to \NN$ is an $\OOO_X/\II$-equivalence if and only if $\NN'/\II\NN' \to \NN/\II\NN$ is an equivalence in $\Mod_{\OOO_X/\II}(\QCoh(X))$, or equivalently in $\QCoh(X)$. Note that when $\II$ is connective, then this can be tested in $\QCoh(V(\II))$, since in this case the vanishing locus $V(\II) \to X$ is affine so that $\QCoh(V(\II)) \simeq \Mod_{\OOO_X/\II}(\QCoh(X))$.
\begin{Def}
\label{Def:I_adic_complete}
    Let $S$ be the class of $\OOO_X/\II$-equivalences. A quasi-coherent $\OOO_X$-module $\KK$ is called \emph{$\II$-complete} if it is $S$-local. Write $\QCoh(X)_{\II}^\wedge$ for the full subcategory of $\QCoh(X)$ spanned by $\II$-complete modules.\footnote{This definition is inspired by \cite[Def.~2.4.10]{HolemanDerivedDelta}, but we follow \cite[\S 7.3]{LurieSpectral} in terminology.}
\end{Def}
\begin{Lem}
    The class $S$ of $\OOO_X/\II$-equivalences is strongly saturated and closed under shifts. Consequently, $\QCoh(X)_{\II}^\wedge$ is presentable and stable, and we obtain a localization
    \[ L_{\II} \colon \QCoh(X) \to \QCoh(X)_{\II}^\wedge \]
    left adjoint to the inclusion $\QCoh(X)_{\II}^\wedge \subset \QCoh(X)$.
\end{Lem}
\begin{proof}
    The first is immediate from the definition, and the second follows from \cite[Prop.~5.5.4.15]{LurieHTT}.
\end{proof}
\begin{Lem}
\label{Lem:internal_map_complete}
    For $\MM \in \QCoh(X)$ and $\KK \in \QCoh(X)_{\II}^\wedge$ it holds that the internal mapping object $\Map(\MM,\KK) \in \QCoh(X)$ is $\II$-complete.  
\end{Lem}
\begin{proof}
    Let $\NN' \to \NN$ be an $\OOO_X/\II$-equivalence. Then also $\NN' \otimes \MM \to \NN \otimes \MM$ is an $\OOO_X/\II$-equivalence. It follows that
    \begin{align*}
        \QCoh(X)(\NN, \Map(\MM,\KK)) & \simeq \QCoh(X)(\NN \otimes \MM, \KK) \\
        & \simeq \QCoh(X)(\NN' \otimes \MM, \KK) \\
        & \simeq \QCoh(X)(\NN', \Map(\MM,\KK)),
    \end{align*}
    whence $\Map(\MM,\KK)$ is $\II$-complete. 
\end{proof}
\begin{Prop}
    The kernel of $L_{\II}$ is a $\otimes$-ideal, so that there is a canonical symmetric monoidal structure on $\QCoh(X)_{\II}^\wedge$ which makes $L_\II$ a symmetric monoidal functor.
\end{Prop}
\begin{proof}
    The first claim is a formal consequence of Lemma~\ref{Lem:internal_map_complete}. The second claim follows from the first by \cite[Thm.~I.3.6]{NikolausTopological}.
\end{proof}

\begin{Lem}
\label{Lem:J_vs_I_adic}
    Let $\JJ \to \II$ be a map of quasi-coherent ideals. Then for all $\MM \in \QCoh(X)$, if $\MM$ is $\II$-complete then it is $\JJ$-complete.
\end{Lem}
\begin{proof}
    This follows from the fact that a morphism $\NN' \to \NN$ in $\QCoh(X)$ which is an $\OOO_X/\JJ$-equivalence is also an $\OOO_X/\II$-equivalence.
\end{proof}
\begin{Lem}[Derived Nakayama]
\label{Lem:I_adically_complete_quotient_zero}
    If $\MM$ is $\II$-complete and $\MM/\II\MM \simeq 0$ then $\MM \simeq 0$.
\end{Lem}
\begin{proof}
    Since $\MM/\II\MM \simeq \MM \otimes \OOO_X/\II$, the assumption implies that $\MM \to 0$ is an $\OOO_X/\II$-equivalence. It follows that
    \[ 0 \simeq \QCoh(0,\MM) \simeq \QCoh(\MM,\MM) \]
    since $\MM$ is $\II$-complete, hence $\MM \simeq 0$.
\end{proof}
\begin{Exm}
    Write $\II_\cl$ for the quasi-coherent ideal on $X$ corresponding to the map $V(\II)_\cl \to X$. By Lemma~\ref{Lem:J_vs_I_adic} it holds that $\II_\cl$-complete implies $\II$-complete. Likewise, if we assume that $V(\II) \to X$ is a closed immersion, then we can make sense of $V(\II)_\red \coloneqq (V(\II)_\cl)_\red$ as the smallest closed subscheme of $X_\cl$ contained in $V(\II)_\cl$, and we write $\sqrt{\II}$ for the ideal corresponding to $V(\II)_\red \to X$. Again $\sqrt{\II}$-complete implies $\II$-complete.
\end{Exm}

\subsection{$\II$-adically complete quasi-coherent $\OOO_X$-modules}
\label{Subsec:Derived_I_cmplt}
Definition~\ref{Def:I_adic_filtration} globalizes as follows. Write 
\[ \DAlg(\Fil^{\geq 0}(-)) \colon \St^\op \to \CAlg(\PrL) \]
for the right Kan extension along the Yoneda embedding of the functor 
\[ \Aff^\op \to \CAlg(\PrL) \colon \Spec(A) \mapsto \DAlg(\Fil^{\geq 0}(\CCC))_{R(A,0)}. \]
Then we obtain an adjunction
\[ (R\dashv \res) \colon \DPair(X) \rightleftarrows \DAlg(\Fil^{\geq 0}(X)) \]
for $X \in \St$. Observe that $R(\OOO_X,\II)$ corresponds to the Rees algebra of a given closed immersion $V(\II) \to X$ in the obvious way via $\Theta$ as classifying stack for filtered objects.

Likewise define $\Fil^{\geq 0}(X) \in \CAlg(\PrL_\st)$ via right Kan extension of $\Spec A \mapsto \Mod_{R(A,0)}(\Fil^{\geq 0})(\CCC)$. Then $\Fil^{\geq 0}(X) \simeq \Fun((\N,\leq,+)^\op, \QCoh(X))$ endowed with Day convolution symmetric monoidal structure and the neutral $t$-structure. By naturality of Kan extensions, we obtain an adjunction
\[ (\ins^n \dashv \ev^n) \colon \QCoh(X) \rightleftarrows \Fil^{\geq 0}(X), \]
where $\ev^n$ is evaluation at $n$ for $n \geq 0$. We also have a monadic adjunction
\[ (\LSym_{R(\OOO_X,0)} \dashv U) \colon \Fil^{\geq 0}(X) \rightleftarrows \DAlg(\Fil^{\geq 0}(X)). \]
\begin{Def}
    We define the \emph{$\II$-adic filtration} of $\MM \in \QCoh(X)$ as
    \[ R(\MM,\II) \coloneqq R(\OOO_X,\II) \otimes_{R(\OOO_X,0)} \ins^0(\MM). \]
    This constitutes a functor
    \[ R(-,\II) \colon \QCoh(X) \to \Mod_{R(\OOO_X,\II)}(\Fil^{\geq 0}(X)). \]
\end{Def}

For a morphism $\II \to \JJ$ in $\Ideal(X)$ we have a map $R(\OOO_X,\II) \to R(\OOO_X,\JJ)$, and thus a map $R(\MM,\II) \to R(\MM,\JJ)$. In particular, we can define $R(\MM,\OOO_X/\II)$ as the cofiber
\[ R(\MM,\II) \to \RR(\MM,\OOO_X) \to R(\MM,\OOO_X/\II) \]
in $\Mod_{R(\OOO_X,\II)}(\Fil^{\geq 0}(X))$. 
\begin{Def}
\label{Def:I_complete}
    The \emph{$\II$-adic completion} of $\MM \in \QCoh(X)$ is 
    \[ \MM^\wedge_{\II} \coloneqq \lim \MM/\II^n \MM \coloneqq \lim R(\MM,\OOO_X/\II), \]
    where the limit is taken in $\QCoh(X)$ after forgetting the $R(\OOO_X,\II)$-module structure. Observe there is a canonical map $\MM \to \MM^\wedge_\II$, natural in $\MM$.\footnote{It is worth mentioning that the smoothness of the derived formalism could be considered to hide the fact that, in the classical case, the formula for $\MM^\wedge_\II$ as quasi-coherent module is in general not well-defined, see \cite[\S 1.3]{BenBerkovich}. The underlying reason this works in the derived setting is that, in this case, the extended Rees algebra is stable under arbitrary base change.}
\end{Def}
\begin{Lem}
\label{Lem:nat_OXI_struct}
    For any $\MM \in \QCoh(X)$ and $n \geq 1$ there is a natural $\OOO_X/\II$-module structure on $\MM/\II^n\MM$ which recovers the $\OOO_X$-module structure upon applying the forgetful functor $\Mod_{\OOO_X/\II}(\QCoh(X)) \to \QCoh(X)$.
\end{Lem}
\begin{proof}
    Since $\MM/\II^n\MM \simeq \MM \otimes \OOO_X/\II^n$ it suffices to do the case $\MM = \OOO_X$, which is clear.
\end{proof}
\begin{Prop}
\label{Prop:completion_is_complete}
    For any $\MM \in \QCoh(X)$ it holds that $\MM^\wedge_\II$ is $\II$-complete.
\end{Prop}
\begin{proof}
    Let $\NN \to \NN'$ be an $\OOO_X/\II$-equivalence. By Lemma~\ref{Lem:nat_OXI_struct} it holds
    \begin{align*}
        \QCoh(X)(\NN',\MM^\wedge_\II) &\simeq \lim \QCoh(X)(\NN', \MM/\II^n\MM)\\
        & \simeq \lim \QCoh(X)(\NN'/\II\NN', \MM/\II^n\MM) \\
        & \simeq \lim \QCoh(X)(\NN/\II\NN, \MM/\II^n\MM)) \\
        & \simeq \QCoh(X)(\NN,\MM^\wedge_\II). \qedhere
    \end{align*}
\end{proof}
\begin{Def}
    A quasi-coherent $\OOO_X$-module $\MM$ is called \emph{$\II$-adically complete} if the canonical map $\MM \to \MM^\wedge_\II$ is an equivalence.
\end{Def}
\begin{Cor}
\label{Cor:derived_complete_implies_adically}
    For any $\MM \in \QCoh(X)$, if $\MM$ is $\II$-adically complete then it is $\II$-complete.
\end{Cor}
\begin{proof}
    Consider the commutative diagram
    \begin{center}
        \begin{tikzcd}
            \MM \arrow[r, "\alpha"] \arrow[d, "\beta"] & \MM^\wedge_\II \arrow[d, "\beta'"] \\
            L_\II(\MM) \arrow[r, "L_\II(\alpha)"] & L_\II(\MM^\wedge_\II).
        \end{tikzcd}
    \end{center}
    Proposition~\ref{Prop:completion_is_complete} shows $\beta'$ is always invertible, hence if $\alpha$ is invertible then so is $L_\II(\alpha)$ and thus $\beta$.
\end{proof}

\begin{Not}
    When $X$ is of the form $X = \Spec(A)$ and $\MM \in \QCoh(X)$ corresponds to $M \in \Mod_A$ and $\II \in\Ideal(X)$ corresponds to $I \to A$, then we write $R(\MM,\II)$ as $R(M,I)$ etc.
\end{Not}

\begin{Exm}
    Consider the setting $\CCC = \Mod_\Z$, and the ideal $\times 0 \colon \Z \to \Z$. Write $(\Z,I_0)$ for the corresponding ideal pair. Then for $M \in \Mod_\Z$ it holds that $R(M,\Z/I_0)$ is the filtered object
    \[ \cdots \xrightarrow{(1,0)} M \oplus M[1] \xrightarrow{(1,0)} M \oplus M[1],\]
    hence that $M \to M^{\wedge}_{I_0}$ is invertible. 
\end{Exm}

\subsection{The case of locally finitely generated ideals}
Recall that a \emph{line bundle} on $X$ is a connective quasi-coherent $\OOO_X$-module which is $\otimes$-invertible. Note that if $V(\II) \to X$ is a virtual Cartier divisor, then (the underlying module of) $\II$ is a line bundle. Indeed, this can be checked on the universal case $B\G_m \to \Theta$, in which case one can show that the associated ideal $\OOO(1)$ corresponds to $\times t^{-1} \colon \mathrm{1}_\CCC[t^{-1}] \to \mathrm{1}_\CCC[t^{-1}]$.\footnote{In the algebraic case $\CCC=\Mod_\Z$ (and with the \'{e}tale topology on $\Aff$), we also have the converse statement. It seems likely this should go through more generally. The missing ingredient is \cite[Prop.~3.2.6]{KhanVirtual}.} 
\begin{Lem}
\label{Lem:VCD_R_lims}
    Suppose that $V(\II) \to X$ is a virtual Cartier divisor. Then the functor $R(-,\OOO_X/\II)$ commutes with limits.
\end{Lem}
\begin{proof}
    From the formula of the extended Rees algebra it follows that $\RR_{V(\II)/X} \simeq \LSym_{\OOO_X}(\II)$,  where the filtration comes about from the maps
    \[ \LSym^{n+1}_{\OOO_X}(\II) \simeq \II^{\otimes (n+1)} \to \II^{\otimes n} \simeq \LSym^n_{\OOO_X}(\II) \]
    obtained by tensoring $\II \to \OOO_X$ with $\II^{\otimes n}$. Then for each $n \geq 0$ the composition
    \[ \QCoh(X) \xrightarrow{R(-,\II)} \Mod_{R(\OOO_X,\II)}(\Fil^{\geq 0}(X)) \xrightarrow{\ev^n} \QCoh(X) \]
    is given by the functor $\MM \mapsto \II^{\otimes n} \otimes \MM$. Since $\II$ is $\otimes$-invertible, the claim follows.
\end{proof}

\begin{Lem}
\label{Lem:iterated_derived_completion}
    Let $\II_1 \to \II$ and $\II_2 \to \II$ be morphisms of connective ideals such that the natural map $V(\II) \to V(\II_1)\times_X V(\II_2)$ is an equivalence, and suppose that $V(\II_2) \to X$ is a virtual Cartier divisor. Then $\MM \to \MM^\wedge_\II$ can be computed as
    \[ \MM \to (\MM_{\II_1}^\wedge) \to (\MM_{\II_1}^\wedge)_{\II_2}^\wedge \simeq \MM^\wedge_\II. \]
\end{Lem}
\begin{proof}
    Put $\BB_i \coloneqq \OOO_X/{\II_i}$ for $i=1,2$ and $\BB \coloneqq \OOO_X/\II$. Since the deformation space commutes with limits, it holds
    \[ \RR^\ext_{\BB/\OOO_X} \simeq \RR^\ext_{\BB_1/\OOO_X} \otimes_{\OOO_X[t^{-1}]} \RR^\ext_{\BB_2/\OOO_X}, \]
    and hence also $\RR_{B_1/\OOO_X} \otimes \RR_{B_2/\OOO_X} \simeq \RR_{\BB/\OOO_X}$ since the functor $(-)^{\ext}$ is symmetric monoidal. We can therefore compute
    \begin{align*}       (\MM_{I_1}^\wedge)_{I_2}^\wedge & \simeq \lim_{k} \ev^k R(\MM_{\II_1}^\wedge,\BB_2) \\
        & \simeq \lim_{k,n} \ev^k R(\ev^n R(\MM,\BB_1),\BB_2) \\
        & \simeq \lim_{i+j \geq 0} \RR^i_{\BB_1/\OOO_X} \otimes \RR^j_{\BB_2/\OOO_X} \otimes \MM \\
        & \simeq \lim R(\MM,\BB),
    \end{align*}
    where we have used Lemma~\ref{Lem:VCD_R_lims} in the second line and a re-indexing argument in the third.
\end{proof}
\begin{Def}
    We say that $\II$ is \emph{locally finitely generated (of rank $n$)} if we can find virtual Cartier divisors $V(\II_j) \to X$ with $1 \leq j \leq n$ and morphisms of ideals $\II_j \to \II$ such that the canonical map
\[ V(\II) \to V(\II_1) \times_X \dots \times_X V(\II_n) \]
is an equivalence. We say that $\II$ is \emph{of the form} $\II = (\II_1,\dots,\II_n)$.
\end{Def}
\begin{Prop}
\label{Prop:derived_I_complete_T}
    Suppose that $\II$ is locally finitely generated of the form $(\II_1,\dots,\II_n)$. Then $\MM \in\QCoh(X)$ is $\II$-adically complete if and only if $\lim R(\MM,\II_j) = 0$ for all $1 \leq j \leq n$.
\end{Prop}
\begin{proof}
    Clearly, $M$ is $\II_j$-adically complete if and only if $\lim R(\MM,\II_j) = 0$. Hence the claim follows from induction on $n$ by Lemma~\ref{Lem:iterated_derived_completion}.
\end{proof}
\begin{Thm}
\label{Thm:derived_I_complete_is_adically}
    Let $\MM \in \QCoh(X)$ be given and suppose that the quasi-coherent ideal $\II$ is locally finitely generated. Then $\MM$ is $\II$-adically complete if and only if it is $\II$-complete.
\end{Thm}
\begin{proof}
    One direction is Corollary~\ref{Cor:derived_complete_implies_adically}. Conversely suppose that $\MM$ is $\II$-complete. Say that $\II$ is of the form $\II= (\II_1,\dots,\II_n)$ for virtual Cartier divisors $V(\II_j) \to X$. Then $\MM$ is also $\II_k$-complete for each $k$, by Lemma~\ref{Lem:J_vs_I_adic}. To show that $\MM$ is $\II$-adically complete, it suffices to show that $\lim R(\MM,\II_k) = 0$ for all $k$. To simplify notation, we may thus assume without loss of generality that $V(\II) \to X$ itself is a virtual Cartier divisor.

    In the exact sequence
    \[ \lim R(\MM,\II) \to \MM \to \MM^\wedge_\II \]
    it holds that $\MM$ and $\MM^\wedge_\II$ are $\II$-complete by Proposition~\ref{Prop:completion_is_complete}, and thus so is $\lim R(\MM,\II)$. By Lemma~\ref{Lem:I_adically_complete_quotient_zero}, it thus suffices to show that $(\lim R(\MM,\II)) \otimes \OOO_X/\II \simeq 0$.
    
    Since $\II$ is $\otimes$-invertible, it holds that $(-) \otimes \II$ commutes with limits, hence so does $(-) \otimes \OOO_X/\II$. We thus obtain 
    \[(\lim R(\MM,\II)) \otimes \OOO_X/\II \simeq \lim R(\MM/\II\MM, \II).\]
    To show that this is zero, we reduce to the case where $X = \Spec A$ and $\II$ is of the form $\times f \colon A \to A$. By shifting degrees, we further reduce to the corresponding statement on $\pi_0$, in which case it is clear.
\end{proof}

\subsection{Formal spectra}
\label{Subsec:formal_spec}
Throughout let still $\II$ be a quasi-coherent ideal on a $\CCC_{\geq 0}$-stack $X$.
\begin{Def}
    The \emph{formal spectrum} of the pair $(\OOO_X,\II)$ is the colimit
    \[ \Spf(\OOO_X,\II) \coloneqq \colim_{n} V(\II^n) \]
    of the $n$-th infinitesimal neighborhoods of the vanishing locus $V(\II) \to X$.
\end{Def}

Let $W \to X$ be a morphism of stacks. Recall that the \emph{pseudo-complement} of $W$ in $X$ is the largest subobject $X \setminus W \hookrightarrow X$ which has empty intersection with $W$.

\begin{Def}
    Write $D(\II) \coloneqq X \setminus V(\II)$. The \emph{formal completion of $X$ along $\lvert V(\II) \rvert$} is the substack 
    \[ X^{\wedge}_{\lvert V(\II) \rvert} \coloneqq X \setminus D(\II) \hookrightarrow X. \]
\end{Def} 

The goal of this subsection is to compare the formal spectrum with the formal completion. We first construct a natural map. Note that there are no assumptions on $\II$ here.

\begin{Lem}
    For all $n \geq 1$ the canonical map $D(\II^n) \to D(\II)$ is invertible. Consequently, we obtain a natural map 
    \[ \Spf(\OOO_X,\II) \to X^{\wedge}_{\lvert V(\II) \rvert} \]
    over $X$.
\end{Lem}

\begin{proof}
    We show the natural monomorphism $D(\II^n) \to D(\II)$ is surjective on $T$-valued points for given $T = \Spec (A) \to X$. Write $I$ for the ideal on $A$ obtained by pulling back $\II$. If $T \to X$ lands in $D(\II)$, then $A/I \simeq 0$, hence $R^{\ext}_{(A/I)/A} \simeq A[t,t^{-1}]$, hence $A/I^n \simeq 0$. The first claim follows.

    The map $V(\II) \to X$ lands in $X^{\wedge}_{\lvert V(\II) \rvert}$, hence also $V(\II^n) \to X$ lands in $X^{\wedge}_{\lvert V(\II^n) \rvert} \simeq X^{\wedge}_{\lvert V(\II) \rvert}$ by the first claim. Passing to the colimit we obtain that $\Spf(\OOO_X,\II) \to X$ indeed factors through $X^{\wedge}_{\lvert V(\II) \rvert} \hookrightarrow X$.
\end{proof}

\begin{Lem}
\label{Lem:Spf_induction}
    Suppose that $\II$ is locally finitely generated of the form $\II= (\II_1,\dots,\II_k)$. Then the natural map
    \[ \Spf(\OOO_X,\II) \to \Spf(\OOO_X,\II_1) \times_{X} \dots \times_X \Spf(\OOO_X,\II_k) \]
    is an equivalence.
\end{Lem}

\begin{proof}
    By Remark~\ref{Rem:inf_nbhd} we obtain a Cartesian diagram
    \begin{center}
        \begin{tikzcd}
            V(\II) \arrow[r] \arrow[d] & V(\II^{n+1}) \arrow[r] \arrow[d] & X \arrow[d] \\
            B\G_m^{\times k} \arrow[r] & (B\G_m ^{\times k})^{(n)} \arrow[r] &  \Theta^{\times k}, 
        \end{tikzcd}
    \end{center}
    where $(B\G_m ^{\times k})^{(n)} \to \Theta^{\times k}$ is the $n$-th infinitesimal neighborhood of the zero section. Hence we reduce to the case where $V(\II) \to X$ is the zero section $B\G_m ^{\times k} \to \Theta^{\times k}$. By construction of quotient stacks, we further reduce to the case where $V(\II) \to X$ is the zero section $\{0\} \to \A^k$. Write $I$ for the corresponding ideal, and put $A \coloneqq \OOO_{\A^k}$. Let $B \in \DAlg_A(\CCC)_{\geq 0}$ be given. We claim that
    \begin{equation}
    \label{Eq:TS_Spf}
        \colim_n \DAlg_A(A/I^n,B) \simeq \colim_n \DAlg_A(A/I_1^n \otimes_A \cdots \otimes_A A/I_k^n ,B),
    \end{equation}
    from which the claim follows by passing to stacks (and a simple cofinality argument).

    The construction of the Rees algebra is natural with respect to morphisms of derived algebraic contexts in the obvious way. In particular, for $F \colon \Mod_\Z \to \CCC$ the unique morphism of derived algebraic contexts it holds that $F(A/I^n) \simeq A/I^n$, where we abusively write also $(A,I)$ for the ideal pair corresponding to $\{0\} \to \A^k$ in the derived algebraic context $\Mod_\Z$. Since $F$ is symmetric monoidal, to show the Equation (\ref{Eq:TS_Spf}) we thus reduce to the case where $\CCC = \Mod_\Z$.

    With the same argument as in Warning~\ref{Warn:dilation_isnt_Rees}, the powers of $I$ are the classical ideal powers, hence $A/I^n$ is the classical quotient ring for all $n$, and likewise for the ideals $I_j$. It follows that $\Spf(\OOO_X,\II)$ is classical, and likewise for $\Spf(\OOO_X,\II_j)$. Moreover, the maps 
    \[ \Spf(\OOO_X,\II_j) \to X \]
    are flat, so that also $\Spf(\OOO_X,\II_1) \times_X \cdots \times_X \Spf(\OOO_X,\II_k)$ is classical.\footnote{We can work with the trivial topology on $\DAlg(\Mod_\Z)_{\geq 0}^\op$ here since we have reduced the statement to the level of prestacks.} Hence to show Equation (\ref{Eq:TS_Spf}) we may assume without loss of generality that $B$ is classical, in which case the statement follows from the familiar cofinality argument.
 \end{proof}

\begin{Rem}
\label{Rem:mild}
    The natural monomorphism 
    \[ [(\A^1 \times \G_m)/\G_m^{\times 2}] \cup [(\G_m \times \A^1)/\G_m^{\times 2}] \to [(\A^2 \setminus \{0\})/\G_m^{\times 2}] \]
    is invertible. This is immediate from the fact that
    \[ [(\A^1 \times \{0\})/\G_m^{\times 2}] \times_{\Theta^{\times 2}} [(\{0\} \times \A^1)/\G_m^{\times 2}]  \simeq B\G_m \times B\G_m, \]
    and the fact that complements of intersections are unions of complements.
\end{Rem}
\begin{Lem}
\label{Lem:Xwedge_induction}
    Suppose that $\II$ is locally finitely generated of the form $\II= (\II_1,\dots,\II_k)$. Then the natural map
    \[ X^{\wedge}_{\lvert V(\II) \rvert} \hookrightarrow X^{\wedge}_{\lvert V(\II_1) \rvert} \cap \dots \cap X^{\wedge}_{\lvert V(\II_k) \rvert} \]
    is an equivalence.
\end{Lem}
\begin{proof}
    By induction we assume $k=2$, and that 
    \begin{align*}
        D(\II_1) &\simeq X \times_{\Theta^{\times 2}} [(\A^1 \times \G_m)/\G_m^{\times 2}] \\
        D(\II_2) &\simeq X \times_{\Theta^{\times 2}} [(\G_m \times \A^1)/\G_m^{\times 2}] \\
        D(\II) &\simeq X \times_{\Theta^{\times 2}} [(\A^2 \setminus \{0\})/\G_m^{\times 2}].
    \end{align*}
    Then by Remark~\ref{Rem:mild} the map $D(\II_1) \cup D(\II_2) \to D(\II)$ is invertible.    
\end{proof}
\begin{Thm}
\label{Thm:Spf_is_completion}
    Suppose that $\II$ is locally finitely generated. Then the natural map
    \[ \Spf(\OOO_X,\II) \to X^{\wedge}_{\lvert V(\II) \rvert} \]
    is an equivalence.
\end{Thm}
\begin{proof}
    By induction on the rank of $\II$, Lemma~\ref{Lem:Spf_induction} and Lemma~\ref{Lem:Xwedge_induction} we reduce to the case where $V(\II) \to X$ is a virtual Cartier divisor. Since the question is local on $X$ we assume that $X$ is affine, say $X = \Spec A$. Write $I$ for the ideal on $A$ corresponding to $\II$. By further localizing, we may assume that $I$ is of the form $\times f \colon A \to A$. Recall that $\DAlg_A(A/(f^n),B)$ is naturally equivalent to the space of paths $f^n \simeq 0$ in $B$. By a similar argument as in the proof of \cite[Lem.~5.1.5]{DAGXII} it follows that the map
    \[ \colim \DAlg(A/(f^n),B) \to \DAlg(A,B) \]
    is an embedding, with as image the full subspace of $\DAlg(A,B)$ spanned by those maps $A \to B$ which send some power of $f \in A$ to $0 \in B$.

    In particular, by passing to stacks we see that $\Spf(\OOO_X,\II) \to X$ is monic as well. Hence it suffices to show that $\Spf(\OOO_X,\II) \to X^\wedge_{\lvert V(\II) \rvert}$ is surjective on $T$-valued points, for $T = \Spec (R) \to X$ given. If $T$ lands in $X^\wedge_{\lvert V(\II) \rvert}$, then $T \times_X D(\II) = \emptyset$, hence the localization $R_f$ vanishes, hence $f^n \simeq 0$ in $R$ for some $n$, hence $T$ factors through $V(\II^n) \to X$ for some $n$. The claim follows.
\end{proof}

\section{Derived $\CCC$-algebra stacks}
We close these notes with a first exploration of transmutation in the setting of geometry relative to derived algebraic contexts. We first investigate stacky filtrations and Rees algebras for animated ring stacks in \S\ref{Subsec:sheaves_of_C_alg} and \S\ref{Subsec:sheaf_of_Rees}. We then introduce transmutations in \S\ref{Subsec:Trans}, give a connection to Weil restriction in \S\ref{Subsec:Weil_CalgSt}, and end with some basic properties of transmutation cohomology in \S\ref{Subsec:trans_coh}.
\subsection{Sheaves valued in adjunctions}
Let $\TTT$ be a site, and consider an adjunction
\[ (F \dashv G) \colon \DDD' \rightleftarrows \DDD \]
of presentable categories.
\begin{Not}
    For a presentable  category $\EEE$, write $\PPP(\TTT;\EEE)$ for the $\EEE$-valued presheaves on $\TTT$, and $\Sh(\TTT;\EEE)$ for the $\EEE$-valued sheaves on $\TTT$.
\end{Not}
Write
\[ (F_*^\PPP \dashv G_*^\PPP) \colon \PPP(\TTT;\DDD') \rightleftarrows \PPP(\TTT;\DDD) \]
for the adjunction induced by postcomposing with $F$ and $G$ respectively.
\begin{Lem}
    Suppose that $F \dashv G$ is monadic. Then also $F_*^\PPP \dashv G_*^\PPP$ is monadic.
\end{Lem}
\begin{proof}
    This follows from Barr--Beck--Lurie since colimits of presheaves are computed section-wise and a split simplicial diagram of presheaves is section-wise split.
\end{proof}
By construction it holds that $G_*^\PPP$ restricts to a functor on sheaves, which we write as $G_*$. By localizing $F_*^\PPP$, we obtain an adjunction
\[ (F_* \dashv G_*) \colon \Sh(\TTT;\DDD') \rightleftarrows \Sh(\TTT;\DDD). \]
\begin{Prop}
\label{Prop:FsharpG_monadic}
    If $F \dashv G$ is monadic, then so is $F_* \dashv G_*$.
\end{Prop}
\begin{proof}
    The corresponding statement on presheaves is immediate from the Barr--Beck--Lurie theorem. Hence it suffices to show that the statement survives to sheaves.

    We have a canonical functor
    \[ \Sh(\TTT;\DDD) \to \PPP(\TTT;\DDD) \times_{\PPP(\TTT;\DDD')} \Sh(\TTT;\DDD') \]
    which we claim is an equivalence. Since it is fully faithful by construction, it suffices to show it is essentially surjective. 

    So let $\varphi \in \PPP(\TTT;\DDD)$ be given such that $G\varphi$ is a sheaf. Then for all $p \in \CCC$ it holds that 
    \[ \CCC(p,G\varphi(-)) \simeq \DDD(Fp, \varphi(-)) \colon \TTT^\op \to \Space \]
    is a sheaf. The bar resolution implies that $\DDD$ is generated under colimits by objects of the form $Fp$. It follows that $\varphi$ is a sheaf, whence the claim.

    Now the functor $G_* \colon \Sh(\TTT;\DDD) \to \Sh(\TTT;\DDD')$ becomes the projection
    \[ \Sh(\TTT;\DDD) \simeq \PPP(\TTT;\DDD) \times_{\PPP(\TTT;\DDD')} \Sh(\TTT;\DDD') \to \Sh(\TTT;\DDD'). \]
    From here, the proposition follows by an abstract argument, using the Barr--Beck--Lurie theorem.
\end{proof}

\begin{Rem}
\label{Rem:triple}
    Suppose that $F \colon \DDD' \to \DDD$ also has a left adjoint $E \colon \DDD \to \DDD'$. Then with the same argument as for $G$ it holds that $F_*^\PPP$ preserves sheaves, hence that $F_*$ is simply the restriction of $F_*^\PPP$. In conclusion, we have adjunctions $E_* \dashv F_* \dashv G_*$ on sheaves, where $F_*$ and $G_*$ are given by postcomposing with $F$ and $G$ respectively, and $E_*$ is induced by postcomposing with $E$ and localization.
\end{Rem}

\subsection{Sheaves of $\CCC$-algebras and their ideals}
\label{Subsec:sheaves_of_C_alg}
Let $F \colon \CCC \to \CCC'$ be a morphism of derived algebraic contexts, with right adjoint $G$. Write also $F \dashv G$ for the adjunction induced on $\DAlg(-)$. Then we obtain a commutative diagram of left adjoint functors
\begin{center}
    \begin{tikzcd}
        \Sh(\TTT;\CCC) \arrow[d, "(\LSym_\CCC)_*"] \arrow[r, "F_*"] & \Sh(\TTT;\CCC') \arrow[d, "(\LSym_{\CCC'})_*"] \\
        \Sh(\TTT;\DAlg(\CCC)) \arrow[r, "F_*"] & \Sh(\TTT;\DAlg(\CCC'))
    \end{tikzcd}
\end{center}
where the vertical functors are the left adjoints of monadic adjunctions by Proposition~\ref{Prop:FsharpG_monadic}.

Applying the above to the morphism $C \colon \Pair(\CCC) \to \Ar(\CCC)$ with right adjoint $K$ from Definition~\ref{Def:ArC} gives us the following. For any morphism $\A \to \B$ in $\Sh(\TTT;\DAlg(\CCC))$ there is a corresponding sheaf of ideal pairs $(\A,\I)$ in $\Sh(\TTT;\DPair(\CCC))$ such that 
the underlying sequence
\[ \I \to \A \to \B \]
is exact in $\Sh(\TTT;\CCC)$. In fact $\B$ can be recovered via $C_*(\A,\I)$ in the sense that
\[ (C_* \dashv K_*) \colon \Sh(\TTT;\DPair(\CCC)) \rightleftarrows \Sh(\TTT;\Arr(\DAlg(\CCC))) \]
is an adjoint equivalence. As before, we write $\A/\I$ for the $\A$-algebra associated to a sheaf of ideal pairs $(\A,\I)$, and $\I(\B)$ for the sheaf of ideals of $\A$ associated to a morphism $\A \to \B$ in $\Sh(\TTT;\DAlg(\CCC))$.

\subsection{The sheaf of Rees algebras}
\label{Subsec:sheaf_of_Rees}
Let $\CCC$ be a derived algebraic context. Applying \S\ref{Subsec:sheaves_of_C_alg} to the constructions in \S\ref{subsec:Rees_Iadic} induces adjunctions
\begin{center}
    \begin{tikzcd}
        \Sh(\TTT;\DAlg(\Fil(\CCC))) \arrow[r, bend right, "(-)^{\geq 0}_*"{name=g}] & \Sh(\TTT;\DAlg(\Fil^{\geq 0}(\CCC))) \arrow[l, bend right, "(-)^{\ext}_*"{name=e}, swap] \arrow[phantom, from=e, to=g, "\dashv"  rotate=-90] \arrow[r, bend right, "\res_*"{name=s}]  & \Sh(\TTT;\DPair(\CCC)), \arrow[l, bend right, "R_*"{name=R}, swap]
        \arrow[phantom, from=R, to=s, "\dashv" rotate=-90]
    \end{tikzcd}
\end{center}
where the functors $(-)_*^{\geq 0}, R_*$ and $\res_*$ can be computed on the level of presheaves---by design for the first one, and by Remark~\ref{Rem:triple} for the latter two.

As before, for a given morphism $\A \to \B$ in $\Sh(\TTT;\DAlg(\CCC))$ with associated ideal $\I$, we define $R_{\B/\A} \coloneqq R_*(\A,\I)$ and $R^\ext_{\B/\A} \coloneqq R_*^\ext(\A,\I)$, where $R_*^\ext \coloneqq (-)^\ext_* \circ R_*$. We consider $R^\ext_{\B/\A}$ a sheaf of derived $\Z$-graded $\CCC$-algebras in the familiar way. Note that $R^\ext_{\A/\A} \simeq \A[t^{-1}]$, where $\A[t^{-1}]$ is the sheaf of derived $\Z$-graded $\CCC$-algebras such that $\A[t^{-1}](X) \coloneqq \A(X)[t^{-1}]$ for all $X \in \TTT$. Hence $R^\ext_{\B/\A}$ is an $\A[t^{-1}]$-algebra in $\Sh(\TTT;\DAlg(\CCC^\Z))$. In particular, we obtain the $\I$-adic filtration
\[ \cdots \to \I^2 \to \I \to \A \]
and derived nilpotent `cothickenings'
\[ \A \to \cdots \to \A/\I^2 \to \A/\I = \B. \]

\begin{Rem}
    Let $\A \to \B$ be a morphism in $\Sh(\Aff(\CCC);\DAlg(\CCC))$ and suppose that the underlying stacks of $\A$ and $\B$ are representable, say $\A = \Spec(A)$ and $\B = \Spec(B)$. Then also the sheaf of ideals $\I \coloneqq \I(\B)$ is representable. Indeed $\I \simeq \Spec (A \otimes_B \mathrm{1}_\CCC)$, where $B \to \mathrm{1}_\CCC$ corresponds to $0 \to \B$. If we write $I$ for the ideal of $B \to \mathrm{1}_\CCC$, then the vanishing ideal of $\I \to \A$ is $I \otimes_BA$.
\end{Rem}

\subsection{Transmutations}
\label{Subsec:Trans}
We now restrict to the case where $\TTT \coloneqq \Aff(\CCC)$ for a derived algebraic context $\CCC$, endowed with a topology as in \S\ref{subsec:geometric_cntxt}. We will describe a version of the transmutation procedure from \cite[Rem.~2.3.8]{BhattF}. To do this, we slightly change perspective. We let $B \in \St(\CCC)$ be given, and suppose we have a functor
\[ \A \colon \Aff(\CCC)_{/B}^\op \to \DAlg(\CCC)_{\geq 0}. \]
In this case, we will always write 
\[ \A \otimes (-) \colon \PPP(\Aff(\CCC))_{/B} \to \St(\CCC) \]
for the unique colimit-preserving extension of the composition
\[ \Spec(-) \circ \A \colon \Aff(\CCC)_{/B} \to \Aff(\CCC) \to \St(\CCC). \] 
Here, $\Aff(\CCC)_{/B}$ is the comma category $\St(\CCC)_{/B} \times_{\St(\CCC)} \Aff(\CCC)$ consisting of $\CCC$-stacks over $B$ which are affine.\footnote{We endow $\Aff(\CCC)_{/B}$ with the Grothendieck topology induced from the one on $\Aff(\CCC)$, so that $\Sh(\Aff(\CCC)_{/B}) \simeq \St(\CCC)_B$.}

\begin{Def}
    The right adjoint of $\A\otimes(-)\colon \PPP(\Aff(\CCC))_{/B} \to \St(\CCC)$ is written
    \[ \St(\CCC) \to \PPP(\Aff(\CCC))_{/B} \colon X \mapsto X^{\A},\]
    and is called \emph{transmutation}.
\end{Def}
Concretely, the transmutation $X^\A$ is the composition 
\[ X^\A \colon \Aff(\CCC)_{/B}^\op \xrightarrow{\A\otimes(-)} \Aff(\CCC)^\op \xrightarrow{X} \Space. \]
\begin{Prop}
\label{Prop:trans_is_stack}
    Let $X \in \St(\CCC)$, and suppose that one of the following conditions holds.
    \begin{enumerate}
        \item $\A$ is an $\DAlg(\CCC)_{\geq 0}$-valued stack and $X$ is affine.
        \item The composition 
        \[ \Spec \A \colon \Aff(\CCC)_{/B}^\op \to \St(\CCC)^\op \]
        is an $\St(\CCC)^\op$-valued stack.
        \item $\A\otimes(-)$ preserves covering families.
    \end{enumerate}
    Then, $X^\A$ is a stack.
\end{Prop}
\begin{proof}
    When $\A$ is an $\DAlg(\CCC)_{\geq 0}$-valued stack, then by definition it holds that
    \[ \DAlg(\CCC)(R, \A(-)) \simeq \Aff(\CCC)(\A\otimes(-),\Spec(R)) \]
    is a stack on $\Aff(\CCC)_{/B}$ for all $R \in \DAlg(\CCC)_{\geq 0}$. In particular this shows (i) when $X = \Spec R$.

    Likewise, when $\Spec \A$ is an $\St(\CCC)^\op$-valued stack, then 
    \[ X^\A(-) \simeq \St(\CCC)(\Spec \A(-), X) \]
    is a stack by definition. This shows (ii).

    Finally, (iii) follows from the adjunction $\A\otimes(-) \dashv (-)^\A$.
\end{proof}

For $f \colon B' \to B$ a morphism of stacks, define the pullback of $\A$ as the composition
\[ f^*\A \colon \Aff(\CCC)_{/B'}^\op \xrightarrow{f_\sharp} \Aff(\CCC)_{/B}^{\op} \xrightarrow{\A} \DAlg(\CCC)_{\geq 0} \]
where $f_\sharp$ is induced by postcomposing with $f$. Since $f_\sharp$ preserves covering families by construction, it holds that $f^*\A$ is an $\DAlg(\CCC)_{\geq 0}$-valued sheaf  if $\A$ is. In fact, $\DAlg(\CCC)_{\geq 0}(R,f^*\A(-)) \simeq f^*\DAlg(\CCC)_{\geq 0}(R,\A(-))$ in $\St(\CCC)_{B'}$ for any $R \in \DAlg(\CCC)_{\geq 0}$. Write also
\[ f^* \A\otimes(-) \colon \Aff(\CCC)_{/B'} \to \Aff(\CCC) \]
for the map corresponding to $f^*\A$.

\begin{Lem}
    For $f \colon B' \to B$ a morphism of stacks and $X \in \St(\CCC)$ it holds that $f^*X^\A \simeq X^{f^*\A}$.
\end{Lem}

\begin{proof}
    By construction and adjunction, for $T \in \Aff(\CCC)_{/B'}$ it holds
    \begin{align*}
        \PPP(\Aff(\CCC))_{/B'}(T, f^*X^\A) & \simeq \PPP(\Aff(\CCC))_{/B}(f_\sharp T, X^\A) \\
    & \simeq  \St(\CCC)(\A\otimes f_\sharp(T),X) \\
    & \simeq \St(\CCC)(f^*\A\otimes T,X) \\
    & \simeq \PPP(\Aff(\CCC))_{/B}(T,X^{f^*\A}). \qedhere
    \end{align*}
\end{proof}

\begin{Lem}
\label{Lem:trans_nat}
    Transmutation is natural in $\A$ in the sense that we have a functor
    \begin{align*}
        \Fun(\Aff(\CCC)_{/B},\DAlg(\CCC)_{\geq 0}) &\to \Fun(\St(\CCC),\PPP(\Aff(\CCC))_{/B}) \\
        \A &\mapsto (X \mapsto X^\A).
    \end{align*}
\end{Lem}

\begin{proof}
    The assignment which sends $\A \colon \Aff(\CCC)_{/B} \to \DAlg(\CCC)_{\geq 0}$ to the associated $\A\otimes(-) \colon \PPP(\Aff(\CCC))_{/B} \to \St(\CCC)$ is clearly functorial in $\A$, and the desired functor follows by passing to right adjoints.
\end{proof}

\subsection{Weil restrictions via $\CCC$-algebra stacks}
\label{Subsec:Weil_CalgSt}
We continue with the case $\TTT \coloneqq \Aff(\CCC)$ for a derived algebraic context $\CCC$. 

Let $V \in \St(\CCC)$ be given such that the projection $p \colon V \to \Spec (\mathrm{1}_\CCC)$ has an affine section $s \colon \Spec (\mathrm{1}_\CCC) \to V$. Consider the functor
\[ \F \colon \Aff(\CCC)_{/V}^\op \to \DAlg(\CCC)_{\geq 0} \colon T \mapsto \Gamma(T \times_V \Spec (\mathrm{1}_\CCC)) \]
which sends a $T$-point of $V$ to the global sections of the pullback of $T$ along $s$. Note that the corresponding functor  $\F\otimes(-) \colon \Aff(\CCC)_{/V} \to \Aff(\CCC)$ preserves covering families so that Proposition~\ref{Prop:trans_is_stack} applies.

For $Y \in \St(\CCC)$ we define the relative version of $\F\otimes(-)$ as
\[ \F_Y\otimes(-) \colon \Aff(\CCC)_{/(Y \times V)}^\op \to \Aff(\CCC)_{/Y}^\op \colon T \mapsto \F\otimes T. \]
For $X \in \St(\CCC)_{/Y}$ we define the relative transmutation as the composition
\[ X^{\F_Y} \colon \Aff(\CCC)_{/(Y \times V)}^\op \xrightarrow{\F_Y\otimes(-)} \Aff(\CCC)_{/Y}^{\op} \xrightarrow{X} \Space. \]

\begin{Prop}
\label{Prop:Weil_via_transmutation}
    For $X \to Y$ in $\St(\CCC)$ it holds that $X^{\F_Y}$ is the Weil restriction of $X \to Y$ along $s\colon Y \times \Spec (\mathrm{1}_\CCC)  \to Y \times V $.
\end{Prop}
\begin{proof}
    By adjunction and construction it holds
    \[ \St(\CCC)_{/(Y \times V)}(T,X^{\F_Y}) \simeq \St(\CCC)_{/Y}(\F\otimes T,X) \simeq \St(\CCC)_{/Y}(T \times_{Y \times V} Y, X)  \]
    for any affine $T$, whence the claim.
\end{proof}
The following generalizes the corresponding statement from the algebraic setting, which was shown in \cite[\S 3.4]{Weil}.
\begin{Exm}
\label{Exm:d2nb_via_transm}
    Consider the derived algebraic context $\CCC^\Z$, and apply Proposition~\ref{Prop:Weil_via_transmutation} to the zero section $s \colon \Spec(\mathrm{1}_\CCC) \to \Spec(\mathrm{1}_\CCC[t^{-1}])$. We use the equivalence $\St(\CCC^\Z) \simeq \St(\CCC)_{/B\G_m}$, so that $s$ corresponds to the zero section $B\G_m \to \Theta$. Consider the functor
    \[ \F \colon \Aff(\CCC)_{/\Theta} \to \DAlg(\CCC)_{\geq 0} \colon \Spec(R) \mapsto R_0 \]
    where $\Spec(R_0) \to \Spec (R)$ is the virtual Cartier divisor classified by the map $\Spec(R) \to \Theta$. Then for any morphism of $\CCC_{\geq 0}$-stacks $X \to Y$, it holds that
    \[ \cD_{X/Y} \simeq (X \times B\G_m)^{\F_{Y \times B\G_m}}. \]    
\end{Exm}

\subsection{Transmutation cohomology}
\label{Subsec:trans_coh}
We continue with the data of a base $B \in \St(\CCC)$ and a functor $\A \colon \Aff(\CCC)_{/B}^\op \to \DAlg(\CCC)_{\geq 0}$. In this subsection we suppose that  $X^\A$ is stack for any $X \in \St(\CCC)$ for simplicity.\footnote{If not then one can sheafify without changing the story, by descent for $\QCoh(-)$.} 
\begin{Def}
    For a $\CCC_{\geq 0}$-stack $X$ and a quasi-coherent module $\EEE \in \QCoh(X^\A)$, we define the \emph{$\A$-cohomology with coefficients in $\EEE$} as 
    \[ H_\A(X,\EEE) \coloneqq \Gamma(X^\A,\EEE) \triangleq p_*(\EEE) \in \QCoh(B) \]
    where $p \colon X^\A \to B$ is the projection map. For $n \in \Z$ we put
    \[ H^n_\A(X,\EEE) \coloneqq \pi_{-n}H_\A(X,\EEE) \simeq \pi_0H_\A(X,\EEE[n]) \in \QCoh(B)^\heartsuit. \]
\end{Def}
Let $f \colon X \to Y$ be a morphism of stacks, write $f^\A \colon X^\A \to Y^\A$ for the associated map. For $\EEE \in \QCoh(Y^{\A})$ put $f^*\EEE \coloneqq (f^\A)^*\EEE$. Then there is a canonical map
\[ H_\A(Y,\EEE) \to H_\A(X,f^*\EEE). \]
We also have the following general descent statement.
\begin{Prop}
    Suppose that $X \simeq \colim X_\alpha$ in $\St$. Let $\EEE \in \QCoh(X^\A)$ and put $\EEE_\alpha \coloneqq \EEE|_{X_\alpha^{\A}}$. Then the natural map
    \[ H_\A(X,\EEE) \to \lim H_\A(X_\alpha,\EEE_\alpha) \]
    is invertible.
\end{Prop}
\begin{proof}
    Write $p \colon X^\A \to B$ and $p_\alpha \colon X_\alpha^\A \to B$ for the structure maps.
    Note that $X^\A \simeq \colim X_\alpha^\A$ on the level of prestacks. By descent for $\QCoh$, it follows that
    \[ \QCoh(X^\A) \simeq \lim \QCoh(X_\alpha^\A). \]
     Under this equivalence, it holds that $p_*$ is the map
    \begin{align*}
        \lim \QCoh(X_\alpha^\A) & \to \QCoh(B) \\
        \{\FFF_\alpha\}_\alpha & \mapsto \lim {p_\alpha}_*\FFF_\alpha,  
    \end{align*}
    whence the claim.
\end{proof}

\begin{Exm}
\label{Ex:ACoh_TCoh}
    Suppose we are in the derived algebraic context $\CCC = \Mod_\Z$, and that $B$ is affine. Then for $\EEE \in \Perf(X^\A)$ and $p \colon X^\A \to B$ the projection it holds
    \begin{equation}
    \label{Eq:Gamma_VB}
    \begin{aligned}
        \tau_{\geq 0} p_*(\EEE) & \simeq \QCoh(B)(\OOO_B,p_*\EEE) \\
        & \simeq \QCoh(X^\A)(\EEE^\vee,\OOO_{X^\A}) \\
        & \simeq \St_{X^\A}(X^\A,\V(\EEE^\vee)).
    \end{aligned}
    \end{equation}
    We can thus understand $\A$-cohomology with coefficients in $\EEE$ as topos-theoretic cohomology, as we now explain.

    Recall that for objects $Y,G$ in a topos $\XXX$, the cohomology of $Y$ with coefficients in $G$ is defined as
    \[ H^0(Y,G) \coloneqq \pi_0\XXX(Y,G). \]
    Further, one says that $G$ has a \emph{delooping} if there is a Cartesian diagram
    \begin{center}
        \begin{tikzcd}
            G \arrow[r] \arrow[d] & * \arrow[d] \\
            * \arrow[r] & BG,
        \end{tikzcd}
    \end{center}
    in which case one puts $H^1(Y,G) \coloneqq H^0(Y,BG)$. This definition is continued inductively to $G$ which have $n$-fold deloopings in the obvious way.\footnote{See, e.g., \cite{NikolausPrincipal}.}

    In particular, in the topos $\XXX = \St_{X^\A}$ we have
    \[       H^0_\A(X,\EEE) \simeq \pi_0 \St_{X^\A}(X^\A,\V(\EEE^\vee)) 
         \triangleq H^0(X^\A,\V(\EEE^\vee)),  \]
    by (\ref{Eq:Gamma_VB}). The $n$-fold deloopings
    $B^n\V(\EEE^\vee) \simeq \V(\EEE^\vee[-n])$
    give
    \[       H^n_\A(X,\EEE) \simeq  H^0_\A(X,\EEE[n]) 
         \simeq H^0(X^\A,\V(\EEE^\vee[-n])) 
         \triangleq H^n(X^\A,\V(\EEE^\vee))    \]
    for all $n \geq 0$.\footnote{Of course, there is nothing special about $X^\A$ here: the argument works for any dualizable $\PPP \in \QCoh(Y)$ on any stack $Y$.}
\end{Exm}

\begin{Con}
\label{Con:transmutation_cofiltration}
    Now suppose we have a morphism $\A \to \B$ with corresponding ideal $\I \to \A$ in $\Fun(\Aff(\CCC)_{/B}^\op,\DAlg(\CCC)_{\geq 0})$. Write $\B_{(k)} \coloneqq \A/\I^k$ for the corresponding nilpotent cothickenings. Let $\EEE \in \QCoh(X^\B)$ be given. Again assume that $X^{\B_{(k)}}$ is a stack for all $k \geq 1$ for simplicity. 

By naturality of transmutation (Lemma~\ref{Lem:trans_nat}) we obtain a factorization
\[ X^\A \to \cdots \xrightarrow{f_3} X^{\B_{(2)}} \xrightarrow{f_2} X^{\B_{(1)}} = X^{\B}, \]
of the natural map $ X^\A \to X^\B$. Let $\EEE \in\QCoh(X^\B)$, and write
\begin{align*}
    \EEE_k \coloneqq \EEE|_{X^{\B_{(k)}}} &&
    \EEE_\infty \coloneqq \EEE|_{X^\A}.
\end{align*}
\end{Con}

\begin{Def}
\label{Def:I_adic_filtr_coh}
    The \emph{$\I$-adic filtration} on the $\A$-cohomology with coefficients in $\EEE_\infty$ is the sequence
    \[ H_{\B}(X,\EEE) \xrightarrow{\varphi_2} H_{\B_{(2)}}(X,\EEE_2) \xrightarrow{\varphi_3} \cdots \to H_{\A}(X,\EEE_\infty) \]
    where $\varphi_k$ is obtained via the unit of the adjunction $f_{k-1}^* \dashv {f_{k-1}}_*$.
\end{Def}

\begin{Def}
    We say that $\A$ is \emph{$\I$-adically complete} if the natural map 
    \[ \A \to \lim \A/\I^k \]
    is an equivalence.
\end{Def}

\begin{Lem}
\label{Lem:Iadic_complete_transmute}
    Suppose that $\A$ is $\I$-adically complete. Then the natural map
     $X^\A \to \lim X^{\B_{(k)}}$ 
    is invertible.
\end{Lem}

\begin{proof}
    This is clear on the level of presheaves, which suffices.
\end{proof}

\begin{Prop}
\label{Prop:exhaustive}
    Suppose that
    \begin{enumerate}
        \item we are in the derived algebraic context $\CCC = \Mod_\Z$,
        \item $B$ is affine,
        \item $\A$ is $\I$-adically complete, and
        \item $\EEE \in \QCoh(X^\B)$ is of the form $\EEE = \PPP|_{X^\B}$ for some $\PPP \in \Perf(B)$.
    \end{enumerate}
    Then the $\I$-adic filtration on $H_\A(X,\EEE_\infty)$ is exhaustive in the sense that the natural map
    \[ \colim H_{\B_{(k)}}(X,\EEE_k) \to H_\A(X,\EEE_\infty) \]
    is invertible.
\end{Prop}

\begin{proof}
    By Lemma~\ref{Lem:Iadic_complete_transmute} it holds $X^\A \simeq \lim X^{\B_{(k)}}$, hence by Example~\ref{Ex:ACoh_TCoh} it holds
    \begin{align*}
        \tau_{\geq 0}H_\A(X,\EEE_\infty) &\simeq \St_{X^\A}(X^\A,\V(\EEE_\infty^\vee)) \\
        & \simeq \St_B(X^\A,\V(\PPP^\vee)) \\
        & \simeq \colim \St_B(X^{\B_{(k)}},\V(\PPP^\vee)) \\
        & \simeq \colim \tau_{\geq 0} H_{\B_{(k)}}(X,\EEE_k),
    \end{align*}
    where we have used that $\V(\PPP^\vee) \to B$ is locally hfp since $\PPP$ is perfect. By shifting $\PPP$ we obtain the corresponding statements on $\tau_{\geq -m}$ for all $m \in \N$, from which we conclude.
\end{proof}

\begin{Exm}
    Let us return to the description of $\cD$ via Weil restrictions from Example~\ref{Exm:d2nb_via_transm}, and we retain that notation but now work in the algebraic setting $\CCC = \Mod_\Z$. Recall that $\Theta$ classifies generalized Cartier divisors $\LL \to \OOO_T$ where $\LL$ is a line bundle on $T$. Write $\OOO(1) \to \OOO$ for the universal generalized Cartier divisor. Then we have a morphism $\V_\Theta(\OOO) \to \F$ of derived rings stacks, with corresponding ideal 
    \[ \V_{\Theta}(\OOO(-1)) \to \V_{\Theta}(\OOO) \to \F. \]

    We want to compute the powers of the ideal $\V_{\Theta}(\OOO(-1)) \to \V_{\Theta}(\OOO)$. Since the powers of ideals of the form $\LL \to \OOO_T$ for a line bundle $\LL$ are given by tensor powers, we see these are just
    \[ \V_{\Theta}(\OOO(-1))^n \simeq \V_{\Theta}(\OOO(-n)). \]

    Write $\F_{(n)} \coloneqq \V_\Theta(\OOO)/\V_{\Theta}(\OOO(-1))^n$ for the corresponding quotient. To ease notation we work implicitly over $Y$. Then the transmutation procedure yields factorizations
    \begin{multline*}
        X \times B\G_m \simeq (X \times B\G_m)^{\V_\Theta(\OOO)}  \to  \dots \to \\ (X \times B\G_m)^{\F_{(2)}} \to (X \times B\G_m)^{\F_{(1)}} \simeq \cD_{X/Y}.
    \end{multline*}
\end{Exm}

\begin{Exm}
    Work in the context $\CCC = \Mod_\Z$ for simplicity. Consider the ideal $\times p \colon \G_a \to \G_a$. For $X \in \St$ and $n \geq 0$ write $X^{[n]}$ for the transmutation of $X$ by $\G_a\sslash p^n$.  Construction \ref{Con:transmutation_cofiltration} yields the filtration
    \[ X \simeq X^{\G_a} \to \dots \to X^{[2]} \to X^{[1]} \simeq X^{\G_a/(p)}. \]
    Observe that for $T = \Spec(C)$ it holds
    \[ X^{[n]}(T) \simeq X(\Spec(C\sslash (p^n))), \]
    hence we can compute
    \begin{align*}
        \lim_n X^{[n]}(T) &\simeq \lim \St(\Spec(C\sslash p^n),X) \\
        & \simeq \St(\colim_n V(p^n),X) \triangleq  \St(\Spf(C,p),X) \\
        & \simeq \St(\Spec(C)^{\wedge}_{\lvert V(p) \rvert} ,X),
    \end{align*}
    by Theorem~\ref{Thm:Spf_is_completion}. 
    In particular, if $C$ is $p$-nilpotent, then $\Spec(C)^{\wedge}_{\lvert V(p) \rvert} \simeq T$, so that $\lim_n X^{[n]}(T) \simeq X(T)$.
\end{Exm}

\bibliographystyle{dary}
\bibliography{refs}	

\end{document}